\documentstyle[11pt, sec, psfig]{article}

\newtheorem{theorem}{Theorem}[section]
\newtheorem{lemma}[theorem]{Lemma}
\newtheorem{proposition}[theorem]{Proposition}
\newtheorem{corollary}[theorem]{Corollary}
\newtheorem{remark}[theorem]{Remark}
\newtheorem{definition}[theorem]{Definition}

\newfont{\Bbb}{msbm10 scaled \magstep1}

\newcommand\bC{\hbox{\Bbb C}}

\newcommand\bR{\hbox{\Bbb R}}
\newcommand\bS{\hbox{\Bbb S}}

\newcommand\bZ{\hbox{\Bbb Z}}

\newfont{\es}{eusm10 scaled \magstep1}
\newfont{\ses}{eufm8 scaled \magstep1}
\newfont{\gt}{eufb10 scaled \magstep1}
\newfont{\sg}{eufb8 scaled \magstep1}
\newfont{\goth}{eufb10 scaled \magstep2}

\newcommand{\re}{\hbox{\gt Re}}
\newcommand{\im}{\hbox{\gt Im}}

\newcommand{\gA}{\hbox{\gt A}}
\newcommand{\gB}{\hbox{\gt B}}

\newcommand{\gog}{{\hbox{\gt g}}}

\newcommand{\gR}{\hbox{\gt R}}

\newcommand{\dir}{\hbox{\es D}}

\newcommand{\en}{\hbox{\es E}}
\newcommand{\uen}{\underline{\hbox{\es E}}}
\newcommand{\hes}{\hbox{\es H}}
\newcommand{\hhes}{\tilde{\hbox{\es H}}}

\newcommand{\Lie}{\hbox{\es L}}
\newcommand{\X}{\hbox{\es X}}

\newcommand{\Z}{\hbox{\es Z}}

\newcommand{\can}{\hbox{\es K}}

\def\ra{\rightarrow}

\def\be{\begin{equation}}
\def\ee{\end{equation}}

\def\lan{\langle}
\def\ran{\rangle}

\newcommand{\na}{{\nabla}}

\newcommand{\hg}{\hat{g}}
\newcommand{\bc}{{\bf c}}
\newcommand{\hbc}{\hat{\bf c}}
\newcommand{\hsi}{{\hat{\sigma}}}
\newcommand{\hpsi}{\hat{\psi}}
\newcommand{\hA}{{\hat{A}}}

\newcommand{\hconf}{\widehat{\hbox{\gt C}}}
\newcommand{\hgauge}{\widehat{\hbox{\gt G}}}
\newcommand{\conf}{\hbox{\gt C}}
\newcommand{\gauge}{\hbox{\gt G}}
\newcommand{\hmodu}{\widehat{\hbox{\gt M}}}
\newcommand{\modu}{\hbox{\gt M}}

\newcommand{\co}{{\sf C}}

\newcommand{\hco}{\hat{\sf C}}

\newcommand{\hbS}{\hat{\hbox{\Bbb S}}}

\newcommand{\ii}{{\bf i}}

\newcommand{\dps}{\dot{\psi}}
\newcommand{\dbps}{{\dot{\bar{\psi}}}}

\newcommand{\dta}{\dot{q}}

\newcommand{\da}{\dot{a}}
\newcommand{\ida}{{\bf i}{\dot{a}}}

\newcommand{\si}{\sigma}

\newcommand{\ve}{{\varepsilon}}
\newcommand{\vfi}{{\varphi}}

\def\dt{\frac{d}{dt}}

\newcommand{\pat}{{\partial_\tau}}

\def\pai{\partial_i}

\def\nah{\hat{\nabla}}

\newcommand{\uso}{\underline{so}}

\begin{document}

\title{Eta Invariants of Dirac Operators on  Circle Bundles Over Riemann Surfaces and Virtual  Dimensions of Finite Energy Seiberg-Witten Moduli Spaces}

\author{Liviu I.Nicolaescu}

\date{June 15, 1997}

\maketitle

\addcontentsline{toc}{section}{Introduction}

\begin{center}
{\bf Introduction}
\end{center}

\bigskip

The eta invariant was introduced in mathematics in the celebrated papers [APS1-3] as a correction  term  in an index formula for a non-local, elliptic boundary value problem and since then it has been   subjected to a lot of scrutiny because of its appearance in many branches of mathematics.

 Contrary to the index density of an elliptic operator, the eta invariant is a {\em non-local} object and  this explains why it is so much harder to  compute. 
   Most concrete computations rely on special topologic or geometric features.  
    For example, one can use the Atiyah-Patodi-Singer theorem to compute  the eta invariant of the signature operator 
    because   in this case the eta invariant is 
    a combination of a  topological term (the signature of a $4k$-dimensional manifold
     with boundary) and a  local contribution  (the integral of the $L$-genus). 
     For $S^1$-bundles over Riemann surfaces,  this approach was successfully carried out in
     [Ko] (see also [O] for similar results in the more general case of Seifert
     manifolds). 
     
For  the   Dirac operator  associated  to a spin structure such an approach   is not possible because  the index of  the  Atiyah-Patodi-Singer problem is notoriously dependent upon the metric.   However, if all the manifolds involved have  positive scalar curvature then a Lichnerowicz type argument  allows the computation of the index and thus, in this case,  the computation   of the eta invariant is a local problem.

The first goal of this paper is to compute the eta invariant  of some Dirac    
operators on the total space of a {\em nontrivial}  circle bundle  $N$ over a  
 Riemann surface  $\Sigma$ of genus $\geq 1$.

  We will do   this for product-like metrics on $N$  such that the fibers   are 
  very short. Such metrics have negative scalar curvatures and thus are beyond 
  the reach of the  Lichnerowicz vanishing approach.    Instead,  using  the   
  recent 
  results of Bismut-Cheeger [BC] and  Dai [Dai]  we will  compute the eta invariant
    for the usual Dirac operator  using  its known adiabatic limit (i.e. its limiting 
    value as the geometry of $N$ changes so that the fibers become shorter and shorter).
    To recover the eta invariant (at least for short fibers)  one can  use  known
     variational formul{\ae}     and some very precise spectral information  about the
      Dirac operators determined by  metrics with very short fibers.   This is  
     entirely a local computation since there are no spectral flow contributions.

Once this computation is performed we embark in a related problem. More precisely,
 we will determine the  eta invariant of a very special scalar perturbation  of the Dirac operator. These perturbed Dirac operators (we called them {\em adiabatic Dirac operators})  arose in  [N] where we studied the  adiabatic limits of the Seiberg-Witten equations on circle bundles (see also [MOY]).  We again use  a variational approach. This time however, there is a   spectral flow contribution  which occurs precisely at  one endpoint of the path of Dirac operators connecting the known situation to the one  for which the computations is sought for.  This requires  some ``spectral care''.    
Additionally,  we indicated  in Appendix C how one can determine  this eta invariant  directly, by elementary means not relying on the Bismut-Cheeger-Dai  theory. This argument extends easily to the more general case of Seifert manifolds.
The eta invariant   is an essential ingredient in the computation of the virtual 
 dimension of the moduli space of finite energy solutions of the Seiberg-Witten 
 equations  on a 4-manifold with boundary a disjoint union of $S^1$-bundles over  Riemann
 surfaces.  The determination of this virtual dimension is an important (albeit intermediary) 
 step in any attempt to establish  gluing formul{\ae} for the Seiberg-Witten   
 invariants. This is the second goal of this paper.

  We compute this dimension  via the Atiyah-Patodi-Singer  and the  Seiberg-Witten 
  analogues of the results in [MMR] describing the structure of the finite energy 
  moduli space.  There is an additional difficulty  one has to overcome.  
The operators  describing the 
 deformation complex of  this moduli space   are based  not just on the adiabatic 
 operator alone. They depend on  a very explicit (though complicated) perturbation of the  direct  
  sum (adiabatic  Dirac operator  + odd signature operator) and the final determination 
  of the virtual dimension requires a very refined  perturbation analysis
  involved in a spectral flow computation.     We obtain explicit formul{\ae}
  for the virtual dimensions  for {\em any 4-manifold  bounding disjoint
  unions of circle bundles}. In particular, one can  immediately determine 
   the virtual dimension of the space of  Seiberg-Witten ``tunnelings''. 
  These are finite energy solutions of the Seiberg-Witten equations  on  an infinite cylinder ${\bR}\times$ circle bundle over a Riemann surface.  This dimension  was also  determined in [MOY] (in the more general context of Seifert manifolds) by identifying the tunnelings with algebraic-geometric objects  on the      orbifold ruled surface associated to a Seifert manifold and then using  a    sophisticated Riemann-Roch argument.  The algebraic-geometric techniques are however ineffective in our more general context and provide no information on the eta invariants which are needed in the computation of Froyshov invariants.    Our methods  extend to the case Seifert manifolds as well but we will pursue this aspect elsewhere.\footnote{In  the paper {\sf dg-ga 9711006} we extended these computations  to Seifert fibrations and used them to  determine the Froyshov invariants of many Brieskorn homology spheres.}   

A  strategy similar to ours  was used in  [SS] to compute the eta invariant of Dirac operators  on circle bundles over Riemann surfaces  of  genus $\geq 2$. 
 There are two main differences.  The first difference  comes from the $spin$  
 structure considered in [SS] which extends to the disk bundle bounding our circle
  bundle.  We perform our computations  on Dirac operators associated to $spin^c$
   structures pulled back from the base of our fibration and these, as explained 
   on p.837 of [KS],  have  notable topological properties. For example, the pullback
    of a $spin$ structure from the base does not extend to a  $spin$-structure on the bounding disk bundle,  though  it extends as a
    $spin^c$-structure.  This explains why the adiabatic limit in [SS] is different from ours and 
    shows that  the eta invariants can distinguish  $spin$ structures!!! The second
 difference is in the manner in which  the adiabatic limit is computed.   In 
 [SS], using  the representation theory of $PSL_2({\bR})$ the authors determine 
 explicitly  the adiabatically important part of the spectrum which allows them
 to determine the adiabatic limit of eta itself.  As we mentioned above, we achieve
  this using the results of Bismut, Cheeger and Dai. In Appendix C we  showed  
  that for the adiabatic Dirac operators the eta {\em function} can be  computed
  directly and   ``elementarily'' and can be elegantly described in terms of   
  Riemann's zeta function  and some topological invariants.

This paper is divided in three sections and  four appendices.  The first  section
 is essentially a brief survey  of known facts concerning  the eta invariant:  
 definition, the Atiyah-Patodi-Singer theorem, variational  formul{\ae} 
 and the   spectral  flow.     We  spent a lot of time  describing these things 
   in an organized way  to eliminate any  ambiguity  concerning the various sign
    conventions. 

The second section contains  the main steps  in the computation of the eta invariants
 discussed  above. We begin by describing the   geometric
 background and the various Dirac operators.  Then using the adiabatic results
 of Bismut-Cheeger-Dai we compute in the second part the eta invariant  of the Dirac operator  on a circle bundle with very short fibers (Theorem \ref{th: eta}). This is a local problem since there are no spectral flows.    In the third part   we compute the eta invariant of the {\em adiabatic Dirac operator} - a perturbation of the Dirac operator which arose
 in [N].   This is  achieved in Theorem \ref{th: etaa} via a variational formula and a spectral flow computation.    The eta invariant of this adiabatic
 Dirac operator is independent of the radius of the fibers. In fact, it depends 
 only on the degree of the fibration and thus   it ``detects'' very little of
 its geometry !!!  This agrees with the way this perturbation was obtained, at 
 the end of a ``shrinking''  process when  at lot of  information is lost.     
 The very general  result of Bismut-Cheeger-Dai may obscure the nice  symmetries
responsible for this unusual conclusion and  we have included in Appendix C an 
elementary  derivation of this result. This elementary approach works in the more
 general context of Seifert manifolds and   enables one to describe the   {\em entire eta function} 
 associated to an adiabatic Dirac operator in terms Riemann's zeta  function.

The last part of this section is devoted to extending the previous computations to the Dirac operators coupled with flat line bundles. We use  essentially the same strategy.  However new phenomena arise  during  the computation of some spectral flow  contributions.   
   
   The third section is devoted to applications to Seiberg-Witten theory.   The
   first two subsections    describe the  3- and 4-dimensional Seiberg-Witten
   equations and some basic facts about them established in [MOY] and [N]. The
   third subsection is entirely devoted   to the computation of a spectral
   flow.  This a very delicate job since one  has to worry about eigenvalues 
   changing sign in a nontransversal  manner.  Detecting these eigenvalues
   requires a very refined perturbation analysis.  In the last subsection  we
   compute virtual dimensions of finite energy Seiberg-Witten moduli spaces on 4-manifolds founding circle bundles over Riemann surfaces. When specialized to  the case of tunnelings (i.e. finite energy solutions on infinite cylinders) we
re-obtain some results of [MOY].

\tableofcontents

\section{The eta invariant of a first order elliptic operator}
\setcounter{equation}{0}

\subsection{Definition}
The elliptic selfadjoint  operators on closed compact manifolds behave in many 
respects as common finite dimensional symmetric matrices. The eta invariant    
 extends the notion of signature  from finite dimensional matrices  to  elliptic
  operators. We will denote the trace of  an  infinite dimensional operator (when it  exists) by  ``Tr'' while  ``tr'' is reserved for finite dimensional operators.  We have the following result.

\begin{proposition}{\rm (a) Consider a closed, compact, oriented Riemann manifold $(N,g)$ of dimension $d$, $E \ra N$ a hermitian vector bundle and $A:C^\infty(E)\ra C^\infty(E)$ a first-order selfadjoint elliptic operator. Then 
\be
\eta_A(s)= \frac{1}{\Gamma(\,\frac{s+1}{2}\,)}\int_0^\infty  t^{(s-1)/2} {\rm Tr}\, (A e^{-tA^2}) dt=\sum_{\lambda>0}\frac{\dim V_\lambda-\dim V_{-\lambda}}{\lambda^s}
\label{eq: 1}
\ee
($V_\lambda=\ker (\lambda-A)$) is  well defined for all ${\re}\, s \gg 0$ and  extends to a meromorphic function  on ${\bC}$.  Its poles  are all simple  and can be  located only at $s=(d+1-n)/2$,  $n=0,1,2, \ldots$.

(b) If $d$ is odd, then the residue of $\eta_A(s)$ at $s=0$ is zero so that $s=0$ is a regular
point.}
\end{proposition}

For a proof of this proposition we refer to [APS3].  When $d$ is odd we define the eta invariant of $A$  by
\[
\eta(A):=\eta_A(0).
\]
\begin{remark}{\rm (a) From the definition it follows directly that 
$\eta(-A)=-\eta(A)$ and $\eta(\lambda A)= \eta (A)$, $\forall \lambda >0$.

(b) In [BF] it is shown that if $A$ is an operator of Dirac type then one  can define its eta invariant directly   by setting $s=0$ in (\ref{eq: 1}). In other words, in this case
\[
\eta(A)=\frac{1}{\sqrt{\pi}}\int_0^\infty t^{-1/2}{\rm Tr}\,(Ae^{-tA^2}) dt.
\]
In the sequel, we will reserve the letter $D$  to denote Dirac type operators.}
\end{remark}

\subsection{The Atiyah-Patodi-Singer theorem}
\begin{figure}
\centerline{\psfig{figure=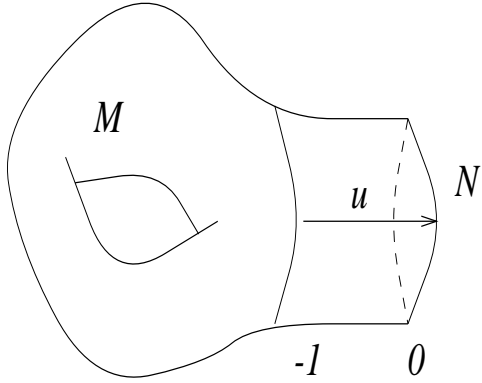,height=2in,width=2.5in}}
\caption{\sl {An oriented manifold with boundary}}
\label{fig: eta1}
\end{figure}

The importance of the eta invariant in mathematics is due mainly to  its appearance  in the formula  for the index of an elliptic boundary value  problem first consider  by  Atiyah-Patodi-Singer in [APS1].

Suppose that $(M,g)$ is  a compact, $(d+1)$-dimensional, oriented  Riemann manifold 
with boundary $N=\partial M$. We assume $d$ is odd  and  that the metric $g$ is a product on a tubular neighborhood $(-1,0]\times N$ of the boundary i.e. $g=du^2+g_0$, where $g_0$ is a  metric on $N$ (see Fig.1).    We orient $N$ such that  the outer normal followed by the orientation of $N$ gives the orientation of
$M$. (This is precisely the orientation that makes the Stokes' formula come out right.)

Next suppose that $E_\pm\ra M$ are two hermitian vector bundles  and $L:C^\infty(E_+)\ra C^\infty(E_-)$ is a first order  elliptic operator which along the neck can be written as
\[
L=G\left(\frac{\partial}{\partial u} - A\right)
\]
where $G:E:=E_+\!\mid_N \ra E_-\!\mid_N$ is a bundle  isomorphism and  $A:C^\infty(E)\ra C^\infty(E)$ is a selfadjoint elliptic operator. (Note that our convention differs from the one in [APS1]!) Denote by
$P_\geq:L^2(E)\ra L^2(E)$ the orthogonal projection onto the closed space spanned by the eigenvectors of $A$ corresponding to eigenvalues
$\geq 0$. $P_<$ is defined similarly. The Atiyah-Patodi-Singer  (APS) boundary value problem is
\[
(APS):\;\;\;\;\left\{
\begin{array}{rcl}
L\psi & = & 0 \\
P_\geq \psi\!\mid_N & = & 0 
\end{array}
\right.
\]
Note that if $\psi$ is a solution of  $(APS)$ then its restriction to the boundary
 lies in the negative eigenspace of $A$. Then,  for all  $u\geq 0$  we
can define
\[
\psi(u) =e^{uA}\psi\!\mid_{\partial M}.
\]
\begin{figure}
\centerline{\psfig{figure=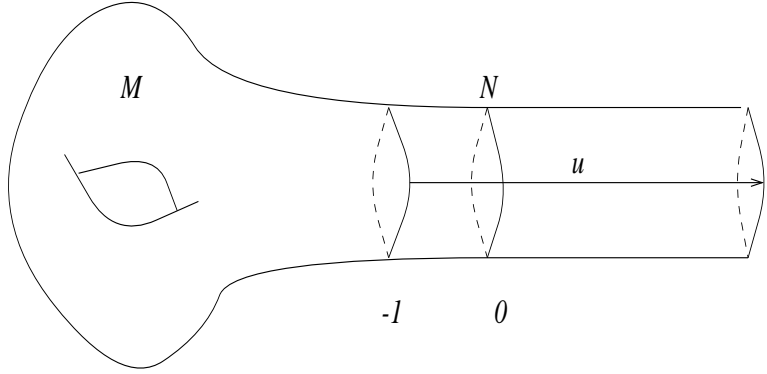,height=2in,width=4in}}
\caption{\sl {Attaching a half-infinite tube}}
\label{fig: eta2}
\end{figure}
We see that  $\psi(u)$ extends $\psi$ to an exponentially decaying solution of $L\psi =0$ on $M_\infty$.  Here $M_\infty$  denotes $M$ with the half-infinite tube $[0,\infty) \times N$ attached (see Fig.2).  Thus, the solutions of $(APS)$  can be identified with the exponentially decaying solutions of $L$ on $M_\infty$. The adjoint  of $(APS)$ is
\[
(APS)^*:\;\;\;\;\left\{
\begin{array}{rcl}
L^*\phi & = & 0 \\
P_< \phi\!\mid_N & = & 0 
\end{array}
\right.
\]
where $L^*$ denotes the formal adjoint  of $L$.  (APS) is an elliptic  problem which implies  finite dimensional spaces of solutions  for both $(APS)$ and its adjoint. Define
\[
{\rm ind}\, (L, APS) =\dim \ker (APS) -\dim \ker (APS)^*.
\]
We have the following fundamental result.

\begin{theorem}{\bf (Atiyah-Patodi-Singer)}
\[
{\rm ind}\, (L, APS)=\int_M\alpha_0(x) dv_g -\frac{1}{2}(\, h(A)+\eta(A)\, )
\]
{\rm where $h(A)=\dim \ker A$,  $\eta(A)$ is the eta invariant of $A$ and $\alpha_0(x) dv_g$ is the index density   determined by $L$ and is  a completely local object (see  [Gky], Sect. 1.8.2  for an exact definition).}
\end{theorem}

Suggested by the above theorem we introduce the $\xi$-{\em invariant}  (or {\em the reduced eta  invariant}) of $A$ by 
\[
\xi(A)= \frac{1}{2}(h(A)+\eta(A)).
\]
Note that $\xi(-A)=(h(A)-\eta(A))/2$ so that  $A\mapsto \xi(A)$ is not an odd function.

In  many  geometrically interesting situations the  index  density $\alpha_0(x) dv_g$ can be described quite explicitly.     We describe below one such instance.

Suppose that $M$ is  equipped with a spin structure. Denote by ${\bS}={\bS}_+\oplus {\bS}_-$ the associated superbundle of spinors.       Fix a connection  $\nabla^M$ on $M$ compatible with  the metric $g$.  $\nabla^M$ {\em need not be the Levi-Civita connection} but we require that it looks like a product in a tubular neighborhood of the boundary.  This  induces in a canonical way  a connection  on ${\bS}$ (compatible with both the metric and the splitting of ${\bS}$) which we  denote by ${\nah}^M$. Suppose moreover that $E\ra M$ is a hermitian vector bundle   equipped  with a  compatible connection $\nabla^E$. We get  in a standard fashion a connection on ${\bS}\otimes E$ compatible with both the metric and the  ${\bZ}_2$-grading.    Finally, this connection canonically defines a Dirac operator  $\hat{\dir}:C^\infty({\bS}_+\otimes E)\ra C^\infty({\bS}_-\otimes E)$ described by
\[
\hat{\dir}: C^\infty({\bS}^+_E)\stackrel{\nabla}{\ra}C^\infty(T^*M\otimes {\bS}_E^+)\stackrel{\hat{\bf c}}{\ra}C^\infty({\bS}_E^-)
\]
where $\hat{\bf c}:T^*M\ra {\rm Hom}\, ({\bS}_+\otimes E, {\bS}_-\otimes E)$ denotes the Clifford multiplication.

As required by the Atiyah-Patodi-Singer index theorem, near the boundary $\hat{\dir}$ has the product structure 
\[
\hat{\dir}=\hat{\bf c}(du)\left( \nabla_u - {\dir}\right)
\]
where ${\dir}$ is the Dirac operator induced by $\hat{\dir}$ on the boundary. 

The index density associated  to this operator is the top degree part of the differential form     $\hat{\bf A}(\nabla^M) \wedge {\bf ch}(\nabla^E)$ where $\hat{\bf  A}$ (resp. {\bf ch}) denote the $\hat{\bf A}$-genus form
(resp.  the Chern character form) obtained from $\nabla^M$ (resp.  $\nabla^E$) via the Chern-Weil construction. In particular, if   $\dim M =4$ and $E$ is the trivial  line bundle equipped with the trivial connection we deduce
\be
{\rm ind}\, (\hat{\dir}, APS)=-\frac{1}{24}\int_M p_1(\nabla^M) -\xi({\dir}).
\label{eq: 2}
\ee

\begin{remark}{\rm The above formula  for $\alpha_0(x) dv_g$ is traditionally proved only  for the
special case when $\nabla^M$ is the Levi-Civita connection. However, a careful
inspection of the proof  in Chap. 11 of [Roe] shows it  extends verbatim to
the more general case when $\nabla^M$ is only metric compatible}
\end{remark}

\subsection{Variational formul{\ae}}
While the eta invariant itself  is a very complex object, its deformation theory turns out to be  a lot simpler. We collect here some results  we will use  in our computations. More specifically, we will address the following problem.   

\medskip

{\it Consider  two metrics $g_i$ $i=0,1$ and compatible connections $\nabla^i$ on 
an odd dimensional manifold  $N$ and denote by ${\dir}_i$ the associated Dirac operators. Compute 
$\xi({\dir}_1)-\xi({\dir}_0)$.}

\medskip

We will soon see this problem  does not have an unique answer and the reason will
 be very  clear.  Leaving this worry aside for a moment, consider a  smooth path
$\{(g_t, \nabla^t)\}_{t\in [0,1]}$ of metrics and compatible connections connecting $(g_0, \nabla^0)$ to   $(g_1, \nabla^1)$. Denote the associated Dirac operators by ${\dir}_t$ and  set
$\xi_t = \xi({\dir}_t)$.  We want to compute  $\dot{\xi}_t =\frac{d\xi_t}{dt}$ 
although at this moment we have
 no guarantee the map $t\mapsto \xi_t$ is differentiable.   

\begin{figure}
\centerline{\psfig{figure=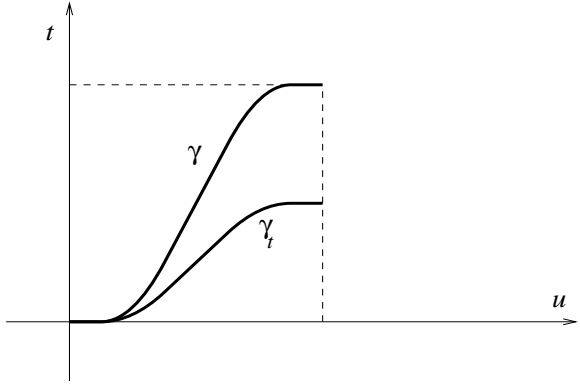,height=2in,width=3in}}
\caption{\sl {The  smoothing function $\gamma$}}
\label{fig: eta3}
\end{figure}

 Since the path $(g_t,\nabla^t)$  may not be  independent of $t$ near $t=0,1$ we
 need to  smooth-out the corners. To this aim, consider a  smooth, nondecreasing map $\gamma: [0,1]\ra [0,1]$, $u\mapsto \gamma(u)$ such that $\gamma(0)=0$, $\gamma(1)$ and $\gamma'(u)\equiv 0$ for $u$ near $0$ and $1$
 (see Figure 3). Moreover for each $0<t\leq 1$ set $\gamma_t(u) =t\gamma(u)$ so that $\gamma_t$ connects $0$ to $t$.

Now for every $0<t\leq 1$ form the operator $L_t$ on $[0,1]\times N$ defined by
\[
L_t =\nabla_u -{\dir}_{t\gamma(u)}.
\]
$L_t$ is an elliptic operator  and  from the A-P-S theorem we get
\[
i_t:={\rm ind}\, (L_t, APS) = \rho_t -\frac{1}{2}(h_0+h_t)+\frac{1}{2}(\eta_0-\eta_t)
\]
where $\rho_t$ denotes the integral of the index density of $L_t$, $h_t=h({\dir}_t)$, $\eta_t=\eta({\dir}_t)$. The above formula can be rewritten as
\be
\xi_t-\xi_0= \rho_t +j_t
\label{eq: var}
\ee
where $j_t=-(h_0+i_t)$. The term $\rho_t$  depends smoothly on $t$ since the coefficients of $L_t$ do. 
 The term  $j_t$ is  ${\bZ}$-valued  so it cannot be smooth, unless it is constant. 
 If $[\xi_t]=\xi_t\; ({\rm mod}\; {\bZ})$ then the map $t \mapsto [\xi_t]$ is smooth
 and by (\ref{eq: var})
 \be
 \frac{d[\xi_t]}{dt} =\dot{\rho}_t.
 \label{eq: var0}
 \ee
   We will deal with $\dot{\rho}_t$ a bit later later but first we need to better understand
   the special nature of the discontinuities of $\xi_t$.

We see from (\ref{eq: 1}) that the  discontinuities of $\xi_t$ (and hence those of $j_t$) are due to jumps in $h_t$.   We describe how the jumps in $h_t$ 
affect $\xi_t$ in a simple, yet generic situation. We assume ${\dir}_t$ is a regular family i.e.

\medskip

\noindent ${\bullet}$ The resonance set  ${\Z}=\{t \in [0,1]\; ;\;  h_t \neq 0\}$ is
finite.

\noindent  ${\bullet}$ For every $t_0\in {\Z}$ there exists ${\ve}>0$, an open neighborhood
 ${\cal N}$ of $t_0$ in $[0,1]$ and smooth maps  $\lambda_k:{\cal N}\ra (-{\ve}, {\ve})$,
 $k=0,1,...,h_{t_0}$ 
 such that for all $t\in {\cal N}$  the  family  $\{\lambda_k(t)\}_k$  describes  {\em all} 
 the eigenvalues of ${\dir}_t$  in $(-{\ve}, {\ve})$ (including multiplicities) 
 and moreover, $\dot{\lambda}_k(t_0)\neq 0$ for all $k=1,2,...,h_{t_0}$.

Now for each  $t\in {\Z}$ set
\[
\sigma_\pm (t) =\#\{ k\;;\; \pm \dot{\lambda}_k(t) >0\}.
\] and
\[
\Delta_t\sigma=\left\{
\begin{array}{rll}
-\sigma_-(0) & ,& t=0\\
\sigma_+(t)-\sigma_-(t)&,& t\in (0,1) \\
\sigma_+(1)&,& t=1 
\end{array}
\right. .
\]
If 
\[
\Delta_t\xi :=\lim_{{\ve}\ra 0^+}(\xi_{t+{\ve}}-\xi_{t-{\ve}})
\]
we see that $\Delta _t\xi=0$ if $t\not\in {\Z}$ while for $t\in {\Z}$ we have 
\begin{figure}
\centerline{\psfig{figure=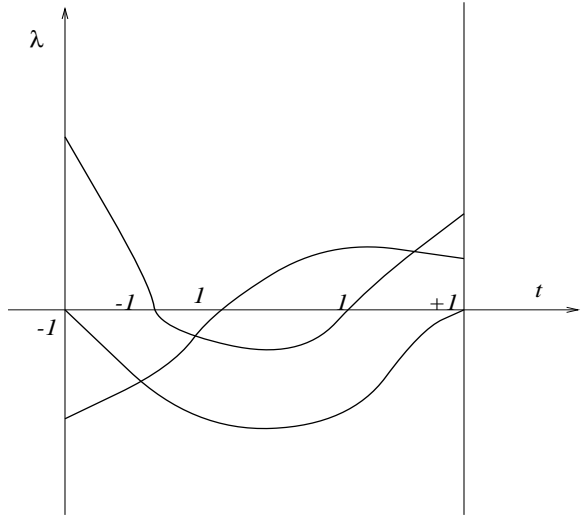,height=2.7in,width=3in}}
\caption{\sl {Spectral flow}}
\label{fig: eta4}
\end{figure}
\be
\Delta_t\xi =\Delta_t\sigma.
\label{eq: sf1}
\ee
(To understand  the above formula  it is convenient to treat ${\dir}_t$ as a finite dimensional  symmetric matrix and then  keep track of the changes in its signature as   the spectrum changes in the
regular way described above).  Finally, define the spectral flow of the family ${\dir}_t$ by
\be
SF({\dir}_t)=\sum_{t\in[0,1]}\Delta_t\sigma.
\label{eq: sf2}
\ee
For example,  in Figure 4    we have  represented those eigenvalues  $\lambda_t$ of
a smooth  path of Dirac operators which  vanish for some values of $t$. The
$\pm 1$'s  describe the jumps $\Delta_t\sigma$. Thus the spectral flow in Figure 4 is 
1.

 Using the equalities $j_1-j_0=\sum_t\Delta_t\xi$   and $j_0=0$ we  deduce
 \be
 j_1-j_0=-i_1-h_0=\sum_t\Delta_t\xi= \sum_{t\in[0,1]}\Delta_t\sigma = SF({\dir}_t) 
 \label{eq: varia}
 \ee
 so that
\be
i_1={\rm ind}\, (L_1, APS)= -h_0 -SF({\dir}_t).
\label{eq: ind}
\ee
From  the equalities  (\ref{eq: var}) and (\ref{eq: varia}) we now conclude
\be
\xi_1-\xi_0= SF({\dir}_t)+\int_0^1 \frac{d[\xi_t]}{dt} dt.
\label{eq: var1}
\ee
\begin{remark}{\em In the above two equalities we  have  neglected the
smoothing effect of $\gamma$. However, since $\gamma(u)$ is nondecreasing, the
crossing patterns of the eigenvalues of $t\mapsto {\dir}_t$ and $u\mapsto
{\dir}_{\gamma(u)}$ are identical.  This implies
$SF({\dir}_t)=SF({\dir}_{\gamma(u)})$.}
\end{remark}

 It is now the time to explain the continuous variation $\frac{d}{dt}[\xi]_t$.
 Formula (\ref{eq: var0}) shows this is a locally computable quantity.  In
 fact,  one can be more accurate than this.   We start with a simple situation first. 
 
 Assume $(N,g)$ is an oriented Riemann
 manifold of dimension $d\equiv 3\; {\rm mod}\; 4$ equipped with a spin structure. Fix a smooth path
 $(\nabla^t)_{t\in[0,1]}$ of $g$-compatible connections and for each $t$ denote
 by ${\dir}_t$ the associated Dirac operator. Consider now the manifold
 $M=[0,1]\times N$ equipped with the metric $\hat{g}=du^2 +g$. The connection
 ${\nah}=du\wedge
 \partial_u +\nabla^{\gamma(u)}$ is compatible with the metric $\hat{g}$ and it
 determines a Dirac operator $\hat{\dir}$ which has the form
 \[
 \hat{\dir}=\hat{\bf c}(du)\left(\partial_u -{\dir}_{\gamma(u)}\right).
 \]
 The A-P-S theorem then gives
 \[
 \xi_1-\xi_0 =\int_M \hat{\bf A}({\nah})-{\rm ind}\,(\hat{\dir},APS) -h_0
 \]
 \[
 =SF({\dir}_{\gamma(u)}) + \int_M \hat{\bf A}({\nah}).
 \]
 To further simplify this formula  note firstly that
 \[
 SF({\dir}_{\gamma(u)})=SF({\dir}_t;\; 0\leq t \leq 1).
 \]
 Secondly, as in [APS2] one can show that the integral term  is {\em independent} of the path
 of connections
 chosen to deform $\nabla^0$ to $\nabla^1$. Thus we can set
 $\nabla^t=\nabla^0+t(\nabla^1-\nabla^0)$.  The resulting integral over $M$ can
 then be rephrased as an integral  over $N$  of the transgression form from
 $\nabla^0$ to $\nabla^1$. This is defined  as the  degree $d$ part of
 \[
 T\hat{\bf A}(\nabla^1, \nabla^0):= \frac{d+1}{2}\cdot  \int_0^1 \hat{\bf A}(\omega,
 \Omega_t) dt
 \]
 where $\omega =\nabla^1-\nabla^0$ and $\Omega_t$ is the curvature of
 $\nabla^0+t\omega$.  More explicitly,
 \[
 \Omega_t=\Omega_0 + td^{\nabla^0}\omega +t^2\omega\wedge \omega
 \]
 where $d^{\nabla^0}$ denotes the exterior derivative defined by $\nabla^0$.
 
  In the special case  when $d=3$ the only important part of
 $\hat{\bf A}$ is $-\frac{1}{24} p_1$ where  $p_1$ is the  degree 2 invariant
 polynomial on $\uso(4)$ given by
 \[
 p_1(X,Y)= -\frac{1}{8\pi^2} {\rm tr}\,(XY).
 \]
 (Here we use the conventions of [BGV]). In this case  the transgression is a multiple of the Chern-Simons integrand,
 and more precisely
 \[
 T\hat{\bf A}(\nabla^1,\nabla^0)= \frac{1}{96\pi^2}{\rm tr}\,
 (\omega\wedge\Omega^0+\frac{1}{2}\omega \wedge d^{\nabla^0}\omega
 +\frac{1}{3}\omega\wedge \omega\wedge \omega).
 \]
 Thus when $d=3$ we have the following remarkable formula
 \be
 \xi_1-\xi_0 = SF({\dir}_t) + \frac{1}{96\pi^2}\int_N {\rm tr}\,
 (\omega\wedge\Omega^0+\frac{1}{2}\omega \wedge d^{\nabla^0}\omega
 +\frac{1}{3}\omega\wedge \omega\wedge \omega).
 \label{eq: trans1}
 \ee
 
 Now consider a more complicated problem.   $N$ is again a compact, oriented,
 $d$-dimensional manifold ($d=3\;({\rm mod})\; 4$), but this time we allow the
 metric to vary. Thus, let $(g_t)_{t\in[0,1]}$ be a smooth  path of Riemann
 metrics on $N$ and, for each $t$, denote by $\nabla^t$ the {\em Levi-Civita}
 connection associated to the metric $g_t$. We obtain in this way  a path of
 Dirac operators  $({\dir}_t)_{t\in [0,1}$. We want to compute $\xi_1-\xi_0$
 assuming for simplicity that all  the operators  ${\dir}_t$ are invertible so
 there is no spectral flow.  
 
 Form again the metric $\hat{g}=du^2
 +g_{\gamma(u)}$ on $[0,1]\times N$ and denote by ${\nah}$ its associated
 Levi-Civita connection. We get a Dirac operator $\hat{\dir}$ on $M$. It has
 the form $\hat{\bf c}(du)\left(\partial_u -{\dir}_{\gamma(u)}\right)$ for $u$
 close to $0$ and $1$. Unfortunately, for $u$ away from the endpoints it has 
 the form
 \[
 \hat{\bf c}(du)\left(\partial_u -{\dir}'_{\gamma(u)}\right)
 \]
 where ${\dir}'_{\gamma(u)}={\dir}_{\gamma(u)} +T_u$   and $T_u$ is a certain
 endomorphism   expressible in terms of $\frac{d}{du}g_{\gamma(u)}$.  If the
 operators ${\dir}'_{\gamma(u)}$  where invertible, then their spectral flow
 would be zero and then  $\xi_1 -\xi_0$ would be expressible as  an integral of
 an $\hat{\bf A}$-form. 
 
 Fortunately, there is a simple way to guarantee the above invertibility, 
 relying on the observation that the size of $T_u$ is comparable with the size
 of the $u$-derivative of $g_{\gamma(u)}$.  Consider  a very large  positive
 number $L$ and  form the tube  $M_L=[0,L]\times N$ equipped with the metric
 \[
 g_L=dv^2 +g_{\gamma(v/L)}.
 \]
 In other words, the path $v\mapsto g_{\gamma(v/L)}$ defines a very, very slow
 deformation of $g_0$   to $g_1$. (A physicist would call this an adiabatic
 process.)  In this  case  the $v$ derivatives of $g_{\gamma(v/L)}$ become 
 extremely small  so that the corresponding perturbations $T_v$ become
 negligible and  ${\dir}'_{\gamma(v/L)}$ will be invertible. If $\nabla^L$ denotes
 the Levi-Civita connection of $g_L$ we get
 \[
 \xi_1-\xi_0 =\int_{M_L}\hat{\bf A}(\nabla^L).
 \]
 As remarked in [APS2], the above integral does not change if we replace
 $\nabla^L$ by a linear connection on $[0,L]\times N$, {\em not necessarily   compatible with}
 $g_L$, which interpolates affinely between the Levi-Civita connections of $g_0$ and
 $g_1$.   This   shows that even  in this case we can express the variation
 of $\xi$ as the integral of a transgression form. The only difference this
 time is that the transgression  goes through $GL(d, {\bR})$-connections rather than
 $O(d)$-connections.  This is  no problem since the two groups are
 homotopically equivalent.
 
 We now have  (almost) all the background  necessary to compute  eta invariants
  of Dirac operators. The only missing piece of information is
  Bismut-Cheeger-Dai  result concerning  the adiabatic limits of the  eta
  invariants. We will state the  special case   we need at the  opportune
  moment.
  
  \begin{remark}{\rm The above observations can be used to determine the index
  of  an elliptic problem   on a {\em noncompact} manifold considered in [LM].

  Consider a smooth, non-decreasing function $\beta:{\bR}\ra [0,\infty)$ such that 
  $\beta(u)\equiv 0$ for $u\leq 1/4$ and $\beta (u)\equiv u$ for $u\geq 3/4$.
  
  Using the notations of \S 1.2,  we define   for each $\mu \in {\bR}$   the
  weighted Sobolev spaces $L^{k,2}_\mu(E_\pm)$   as completions of
  $C_0^\infty(E_\pm)$ with respect to the norm
  \[
  \|\psi\|_{k,2,\mu}=\left(\sum_{j=0}^k\int_{M_\infty}|e^{\mu\beta(u)}\nabla^j\psi|^2
   dV_g\right)^{1/2}.
  \]
   Consider the bounded operator $L=\partial_u-A:L^{1,2}_\mu(E_+)\ra
   L^2_\mu(E_-)$. In [LM] it was shown that   $L$ is Fredholm if and only if 
   $A+\mu$ is invertible, i.e. $-\mu \not\in {\rm spec}\,(A)$. We denote by
   $i_\mu(L)$ its index.  For example, if $A$ is invertible, then as pointed out
   in [APS1] we have
   \[
   i_0={\rm ind}\,(L, APS).
   \]
  In general, to compute $i_\mu$ for an arbitrary $\mu$  note that the map
   \[
  T_\mu: L^{2}_\mu \ra L^{2},\;\;\;\psi \mapsto e^{\beta(u)}\psi
  \]
  is an isometry so  that $i_\mu (L)=i_0(T_\mu L T^{-1}_\mu)$.  A simple
  computation shows that
  \[
  T_\mu L T^{-1}_\mu=L_\mu:= L-\mu\beta'(u).
  \]
  Construct $M_1$ by attaching the cylinder $C_1=[0,1]\times N$ to the boundary of
  $M$. Alternatively, $M_1$ is the region $u\leq 1$ in $M_\infty$.  Then
  $L_\mu$ is well defined on $M_1$ and  as above we conclude
  \[
  i_\mu= {\rm ind}\,(L_\mu, APS).
  \]
  Set $A_\mu = A+\mu$.   We have
  \[
  {\rm ind}_{M_1}\, (L_\mu, APS) -{\rm ind}_M\,(L, APS)= -(\,
  \xi(A_\mu)-\xi(A)\,)
 +\int_{C_1}\alpha_0(x)dv_g.
 \]
 On the other hand, the  above index density can be expressed  as in (\ref{eq:
 var}) in terms of the APS index of the operator $L-\mu\beta'(u)$ on $C_1$.
 \[
 \int_{C_1}\alpha_0(x)dv_g=\xi(A_\mu)-\xi(A) +h(A) +  {\rm ind}_{C_1}\,
 (L-\mu\beta'(u), APS)
 \]
 Finally, according to (\ref{eq: ind}), the  last term can be expressed as a spectral flow
 \be
 {\rm ind}_{C_1}\,(L-\mu\beta'(u), APS)=-h(A) -SF(A+t\mu,\; t\in[0,1]).
 \label{eq: ex1}
 \ee
 Putting all of the above together we obtain the following useful equality
 \be
 i_\mu={\rm ind}\,(L, APS) -SF(A+t\mu,\;t\in[0,1])
 \label{eq: lm}
 \ee
 This is in perfect agreement with Theorem  1.2 in [LM]. Note also that if
 $\mu$ is sufficiently small  and positive then there is no spectral flow correction in the
 above formula.
 
 We also want to mention an immediate consequence of the above considerations.
 Consider two elliptic  first order operators  as above
 \[
 L_1,L_2:\Gamma(E_+)\ra \Gamma(E_-)
 \]
 which  have the {\em same principal symbol} and have the $APS$ form $L_j=G(\partial_u-A_j)$
 along the neck. Arguing as in the proof  of (\ref{eq: lm}) we deduce the
 {\em excision formula}
 \be
 {\rm ind}\,(L_2,APS)={\rm ind}\,(L_1, APS) -SF(A_1\ra A_2)
 \label{eq: ex2}
 \ee
 where $SF(A_1\ra A_2)$ denotes the spectral flow of the affine path of
 elliptic operators $A_t=A_1+t(A_2-A_1)$, $t\in[0,1]$.}
 \label{rem: 110}
 \end{remark}

\section{Eta invariants of Dirac operators}
\setcounter{equation}{0}

\subsection{The differential geometric  background}
Consider $\ell \in {\bZ}$ and  denote by $N=N_\ell$ the total space of a degree
$\ell$ principal $S^1$ bundle over a  compact oriented surface of genus $g$: $S^1 \hookrightarrow N_\ell \stackrel{\pi}{\ra} \Sigma$.  Denote by $\zeta\in{\rm Vect}\,(N)$ the infinitesimal generator of the $S^1$ action.  $N$ has a natural orientation  which can be described using any   any splitting $TN=\lan \zeta\ran \oplus \pi^*T\Sigma$ determined by an arbitrary connection.

Assume $\Sigma$ is equipped with a  {\em constant sectional curvature} Riemann metric $h_b$ such that ${\rm vol}_{h_b}(\Sigma)= \pi$. Pick a connection form ${\bf i}\vfi\in {\bf i}\Omega^1(N)$ such that
 \[
 -d\vfi= 2 \ell dv_{h_b}.
 \]
This choice is possible since $\frac{-1}{2\pi}d\vfi$ represents the first Chern class of $N$ which is $\ell$.  For each $0<r \leq  1$ define a metric $h_r$ on $N$ by
\[
h_r=\vfi_r\otimes \vfi_r \oplus \pi^* h_b,\;\; \vfi_r=r\vfi.
\]
Set $\zeta_r=r^{-1}\zeta$.

Using this metric we can orthogonally  split $T^*N\cong \lan \vfi \ran \oplus \pi^*T^*\Sigma$  and  this defines in a natural way an orientation on $N$.  If $\ast_r$ denotes the Hodge $\ast$ operator of the metric $h_r$ we get 
\be
d\vfi_r =2\lambda_r \ast_r \vfi_r
\label{eq: di0}
\ee
where $\lambda_r=-r\ell$.

Fix  a local, orthonormal  coframe  $\theta^1, \theta^2$ on the base $\Sigma$
such that 
\be
d\theta^1=\kappa\theta^1\wedge \theta^2
\label{eq: di1}
\ee
 and 
\be
d\theta^2=0
\label{eq: di2}
\ee
 where $\kappa$
is a {\em nonnegative constant}.  Such a choice is obviously possible if $\Sigma$ is the flat torus. If $\Sigma$ has higher genus then any  constant
curvature  metric on  $\Sigma$ admits such local coframes because it is a
quotient of the hyperbolic plane and on the hyperbolic plane such choices are
possible.   Note that $-\kappa^2$ is actually the sectional curvature of $\Sigma$
so by Gauss-Bonnet
\[
\kappa^2 =4(g-1).
\]
We  now get a local,  oriented, orthonormal frame  of $T^*N$
$(\vfi_r,\vfi^1,\vfi^2)$, where $\vfi^i=\pi^* \theta^i$, $i=1,2$. Denote
 by $(\zeta_r,\zeta^1,\zeta^2)$ its dual frame.  In [N] we showed that the
 1-form associated to the Levi-Civita connection $\nabla^r=\nabla(h_r)$ by the above frame is
 \be
 \omega_r= \left[
\begin{array}{rrr}
0 & A_r & -B_r \\
-A_r & 0 & C_r \\
B_r & -C_r & 0 
\end{array}
\right] 
\label{eq: structural}
\ee
where 
\be
A_r=\lambda_r\vfi^2,\;\;B_r=\lambda_r\vfi^1,\;\;C_r=-\lambda_r
\vfi_r -\kappa \vfi^1.
\label{eq: str1}
\ee

In our  computations  we  will need various other connections compatible with 
the above metric.  For $t\in (0,1]$ consider the bundle isomorphism $L_t:TN\ra
TN$ described locally by
\[
\zeta \mapsto t\zeta,\;\;\zeta_i\mapsto \zeta_i,\;\;i=1,2.
\]
Clearly $L_t$ defines an isometry $(TN, h_{rt}) \ra (TN, h_r)$, for all $r>0$
and all $t\in (0,1]$.  This implies that  the  connection
\[
\nabla^{r,t}=L_t \nabla^{rt} L_t^{-1}
\]
is compatible with the metric $h_r$. A simple computation shows that the 1-form
associated to this connection by the frame $(\zeta_r, \zeta_1, \zeta_2)$ is
\be
\omega_{r,t}=\left[
\begin{array}{ccc}
0 & \lambda_rt\vfi^2 & - \lambda_rt \vfi^1 \\
-\lambda_rt\vfi^2  & 0  & -\lambda_rt^2 \vfi_r -\kappa \vfi^1 \\

 \lambda_rt \vfi^1 &\lambda_rt^2 \vfi_r +\kappa \vfi^1 & 
 \end{array}
 \right]
 \label{eq: str2}
 \ee
 The connection $\nabla^{r,t}$   a connection $\tilde{\nabla}^{r,t}$ on  $\lan
 \zeta_r\ran^\perp=\pi^* T\Sigma$.  Denote by $\nabla^\Sigma$ the Levi-Civita
 connection on $\Sigma$. The formulae (\ref{eq: str1}) and (\ref{eq: str2}) imply immediately
 \[
 \lim_{t\ra 0}\tilde{\nabla}^{r,t}= \pi^*\nabla^\Sigma,\;\;\forall r>0
 \]
 and
 \[
 \lim_{r\ra 0}\tilde{\nabla}^{r,t}=\pi^*\nabla^\Sigma, \forall t>0.
 \]
 Set $\nabla^\infty= \pi^*\nabla^\Sigma$.

 The bundle $\lan \vfi\ran^\perp$ has a natural complex structure  locally defined by the correspondences $\vfi_1\mapsto-\vfi_2 \mapsto-\vfi_1$.  
 In this way we get a complex line bundle ${\can}\ra N$. It is isomorphic with 
 the pullback of the canonical line bundle $K_\Sigma$ of the base.   Once we fix a
 $spin$ structure on $\Sigma$ the Levi-Civita induced connection  defines a    
 natural holomorphic structure on $K_\Sigma^{-1/2}$. In [N] we showed   that   
 this induces a $spin$ structure on $N$ with associated bundle of spinors
 \[
 {\bS}={\can}^{-1/2}\oplus {\can}^{1/2}.
 \]
 The Clifford multiplication is described explicitly in Appendix D.  We only
 want to mention here that our choice is such that the Clifford multiplication
 by the volume form on $N$ is equal to $-1$. This agrees with the conventions
 of [BC].
  
  Note that the $h_r$-compatible connections $\nabla^{r,t}$ define  a 2-parameter
 family of Dirac operators ${\dir}_{r,t}$ on ${\bS}$. To   explicitly describe
 their form introduce as in  [N] the following operators.
 \[
 Z =Z_{r}=\left[
 \begin{array}{rr}
 {\bf i}\partial_{\zeta_r} & 0 \\
 0        & -{\bf i}\partial_{\zeta_r}
 \end{array}
 \right]
 \]
 \[
 T=\left[ 
 \begin{array}{cc}
 0  &  {}^\flat\bar{\partial} \\
 {}^\flat\partial  & 0 
 \end{array}
 \right]
 \]
 where (locally) 
 \[
 {}^\flat\bar{\partial}=2^{-1/2}(\vfi^1-{\bf i}\vfi^2)\otimes
 (\nabla^\infty_{\zeta_1}+{\bf i}\nabla^\infty_{\zeta_2})
 \]
  and
  \[
  {}^\flat{\partial}=2^{-1/2}(\vfi^1+{\bf i}\vfi^2)\otimes
 (\nabla^\infty_{\zeta_1}-{\bf i}\nabla^\infty_{\zeta_2}).
 \]
 The computations in [N]  show that
 \be
 {\dir}_{r,t}=Z_{r} +T +\frac{\lambda_rt^2}{2}.
 \label{eq: dir1}
 \ee
 We set
 \[
 {\dir}_r:={\dir}_{r,t}\!\mid_{t=1}
 \]
  and 
 \be
 D_r= \lim_{t\ra 0}{\dir}_{r,t}= Z_{r}+T ={\dir}_r -\frac{\lambda_r}{2}.
 \label{eq: dir2}
 \ee
 The goal of this paper is to compute $\eta({\dir}_r)$ and $\eta(D_r)$.      

We conclude this subsection by listing  properties of the spectrum of $D_r$ as $r\ra 0$. Their proofs can be found in [N].

\begin{proposition}{\rm   (a)  For all $r\in(0,1]$ 
\[
\dim \ker D_r =2h_{1/2}:=2\dim H^0(\Sigma, K^{1/2}).
\]
(b) There exists $r_0>0$ and $z_0>0$ such that for all $r\in (0,r_0]$ the only eigenvalue of $D_r$ in $(-z_0,z_0)$ is  0.}
\label{prop: spec}
\end{proposition}

\subsection{The eta invariant of the  spin Dirac operator}
 In the sequel we assume $\ell\neq 0$ i.e. $N$ is a {\em nontrivial} circle bundle.

As we mentioned in the introduction,  the key step in our computation of $\eta_r:=\eta({\dir}_r)$ will be the adiabatic result of  Bismut-Cheeger [BC]  in the more accurate form of  [Dai].   Instead of formulating the most general version  of their result (which would require   a large preamble) we  state it for the special case we have in mind.  Fortunately, in this case concrete computations were performed  in [Z] and [DZ].   Set $\eta_0=\lim_{r\ra 0}\eta_r$. We then have the following result.

\begin{theorem}{\rm  The adiabatic limit exists and  moreover
\[
\eta_0 =-2\int_\Sigma\hat{\bf A}(\nabla^\Sigma)\frac{\tanh(c/2)-c/2}{c\tanh(c/2)}
+\sum_{\mu} {\rm sign}\,\mu
\]
where $c\in H^2(\Sigma, {\bR})$ is the  Euler class of  the $S^1$ bundle $N$ and the summation is  carried over all  nonzero eigenvalues
$\mu$ of ${\dir}_r$ which are of size $O(r)$ as $r\ra 0$.}
\label{th: z}
\end {theorem}

Using the Taylor expansion
\[
\tanh x = x-\frac{x^3}{3} +O(x^5)
\]
 and the fact that $\hat{\bf A}(\Sigma)=1$ we deduce
\[
{\bf A}(\nabla^\Sigma)\frac{\tanh(c/2)-c/2}{c\tanh(c/2)} =-c/12.
\]
Since ${\dir}_r=D_r+\frac{\lambda_r}{2}$ we deduce  from Proposition \ref{prop: spec}
that for $0<r \ll r_0$ the only nonzero eigenvalue of ${\dir}_r$ which is of size $O(r)$ is
$\lambda_r/2=-r\ell/2$ and it has multiplicity $\dim \ker D_r$. Putting all the above together we deduce
\be
\eta_0 = \ell/6 -2{\rm sign}\,(\ell ) h_{1/2}.
\label{eq: eta}
\ee
Let $r_0$ as in Proposition \ref{prop: spec}.  For every $0<r \ll r_0$ set $\xi_r=\xi({\dir}_r)$. Note that $\xi_r=\frac{1}{2}\eta_r$ since by Proposition \ref{prop: spec}  $\ker {\dir}_r=0$ if $r \ll r_0$.    Finally, denote $\xi_0=\lim_{r\ra 0} \xi_r$. From the  equality  (\ref{eq: eta}) we deduce
\be
\xi_0-\xi_r= \ell/12 -{\rm sign}\, (\ell)h_{1/2} -\xi_r,\;\;\forall r\ll r_0
\label{eq: xi}
\ee
On the other hand, as explained in subsection 1.3,  for all  $0<\rho <r\ll r_0$ the difference  $\xi_\rho-\xi_r$  can be expressed as the integral of an $\hat{\bf A}$-genus form. There is no spectral
flow because the operators ${\dir}_r$ are invertible  for $ 0<r \ll r_0$. 
Following the prescriptions at the end of  subsection 1.3 one obtains the following result. (For details we refer to Appendix A.)

\begin{lemma}{\bf (First  transgression formula)}{\rm  For all $ 0 <r \ll r_0$ we have
\[
\lim_{\rho\ra 0}(\xi_\rho-\xi_r) =- \frac{\ell}{12}(\ell^2r^4-\chi r^2)
\]
where $\chi=\chi(\Sigma)=2-2g$.}
\end{lemma}

By combining all of the above we get the following result.

\begin{theorem}{\rm For all  $0<r\ll r_0$ we have}
\[
\frac{1}{2}\eta_r=\xi_r= \frac{\ell}{12}-{\rm sign}\,(\ell)h_{1/2} +\frac{\ell}{12}(\ell^2r^4-\chi r^2).
\]
\label{th: eta}
\end{theorem}

\subsection{The eta invariant of the adiabatic Dirac operator}
In this  subsection we take  up the computation of the eta invariant of
$D_r$.  We rely on our freshly acquired knowledge of $\xi({\dir}_r)$. 

The Dirac operator ${\dir}_r$ is associated to the Levi-Civita connection $\nabla^r$
 while $D_r$ is associated to the connection  $\nabla^{r,0}=\lim_{t\ra 0}
 \nabla^{r,t}$.     Set $\tau = (1-t^2)$. Then, using  (\ref{eq: dir1}) and
 (\ref{eq: dir2}) we get a path
 \[
 \tilde{\dir}_{r,\tau}= {\dir}_r -\frac{\tau\lambda_r}{2}
 \]
 such that $\tilde{\dir}_{r,0}={\dir}_r$ and $\tilde{\dir}_{r,1}=D_r$. Set
 $\xi_\tau=\xi(\tilde{\dir}_{r,\tau})$. Using the variational
 technique described in \S 1.3 we deduce
 \be
 \xi_1=\xi_0 + SF(\tilde{\dir}_{r,\tau};\, \tau\in[0,1]) +\int_N
 TA(\nabla^{r,0}, \nabla^{r,1}).
 \label{eq: vara}
 \ee
 To compute the spectral flow note that according to Proposition \ref{prop:
 spec} the operator $\tilde{\dir}_{r,\tau}$ has a kernel only for $\tau=1$. In
 this case, the kernel has dimension $2h_{1/2}$.
 Using (\ref{eq: sf1}) of \S 1.3 we deduce
 \be
 SF=\left\{
 \begin{array}{rll}
 2h_{1/2} &, & \ell >0 \\
 0 &, & \ell <0
 \end{array}
 \right.
 \label{eq: sfd}
 \ee
 
 As for the transgression term, it is described in the following lemma whose
 proof can be found in Appendix B.
 
 \begin{lemma}{\bf (Second transgression formula)} 
 \[
 \int_N
 TA(\nabla^{r,0}, \nabla^{r,1}) =-\frac{\ell}{12}(\ell^2r^4-\chi r^2).
 \]
 \end{lemma}
 
 Putting together all the above we obtain  from Theorem \ref{th: eta} the following result.
 
 \begin{theorem}{\rm   For all $0<r \ll 1$ we have}
 \[
 \xi(D_r)=\frac{\ell}{12}+ h_{1/2}.
 \]
 \label{th: etaa}
 \end{theorem}
 
 Note  that $\xi(D_r)$ is independent of $r$ !!! In hindsight, this should not
 be so surprising   if we think that $D_r$ was obtained  after the adiabatic
 deformation in (\ref{eq: dir2}). Notice that $\xi(D_r)$ still ``remembers'' it
 came from a fibration  due to the term $\ell/12$. The geometry of the base is
 reflected in the term $h_{1/2}$. Remarkably, $\eta(D_r)=\ell/6$.  Thus  the 
 base $\Sigma$ is ``invisible'' to the eta invariant of $D_r$ !!!

\subsection{The eta invariant of the coupled adiabatic Dirac operator}
Let us begin by recalling that the Gysin exact sequence implies that
\[
H^2(N;{\bZ})\cong \pi^*H^2(\Sigma;{\bZ}) \oplus H^1(\Sigma,;{\bZ})\cong {\bZ}_{|\ell|}\oplus {\bZ}^{2g}.
\]
Consider  a complex line bundle  $L\ra N$ such that $c_1(L)=\hat{k}\in
{\bZ}_{|\ell|}$. Such a line bundle can be obtain as the pullback of line bundle $L_\Sigma \ra
\Sigma$ of degree $k\in {\bZ}$.   Note that $k$ is determined only modulo
$\ell$.  A line bundle as above admits flat connections and the holonomy of such a connection is $\exp(2\pi k{\bf i}/\ell)$. The collection of gauge equivalence classes of flat connections is homeomorphic to a torus $T^{2g}$. 

These facts were proven in [N]  relying on a simple observation  which we repeat here  since it is relevant to our computations.

Let $A$ be a flat connection on $L$ and set 
\be
B:=A+\frac{k{\bf i}}{\ell} \vfi.
\label{eq: pull}
\ee  
Then $B$ is a connection with trivial holonomy along fibers  and it can be regarded  as a pullback of a connection on a line bundle $L_\Sigma \ra
\Sigma$ such that $c_1(L_\Sigma)=k\in {\bZ}$. Now set
\[
{\bS}_L={\bS}\otimes L ={\can}^{-1/2}\otimes L \, \oplus \, {\can}^{1/2}\otimes L.
\]
By coupling the connection $\pi^*\nabla^\Sigma$ on ${\bS}$  with the flat connection $A$ we get a connection on ${\bS}_L$ which leads to a Dirac operator $D_{A}=D_{A,r}$. We call this the  {\em adiabatic operator coupled with} $A$.

Similarly, using the connection $B$ on $L$ we obtain two connections on ${\bS}_L$ obtained by coupling $B$ with the Levi-Civita connection and respectively the connection $\pi^*\nabla^\Sigma$. These lead to two Dirac operators, ${\dir}_{B,r}$ and  respectively $D_{B,r}$. The goal of this section is to compute the eta invariant of the operator $D_{A,r}$  which  played a key role in [N] in the description of the reducible adiabatic solutions of the Seiberg-Witten equations.  We will use  these eta invariant computations in a forthcoming work on Seiberg-Witten equations on manifolds with cylindrical ends.

The computation of $\eta(D_{A,r})$ for $ 0< r \ll 1$ is performed in three steps.

\noindent {\bf Step 1:}  Compute $\xi({\dir}_{B,r})$.

\noindent {\bf Step 2:} Compute  $\xi(D_{B,r})$.

\noindent {\bf Step 3:} Compute $\xi(D_{A,r})$.

While  the first two steps   follow closely \S 2.2 and  respectively \S 2.3,  
interesting new  phenomena arise  at Step 3.

Before we carry out the computations we need  to review some facts and
introduce  some notations.

Recall first that if $L_\Sigma\ra \Sigma$ is a complex line bundle then any connection
$B$ on $L_\Sigma$ introduces a holomorphic structure on $L_\Sigma$.  We denote by
$h_{1/2}(L_\Sigma)$ the dimension of the space of holomorphic sections  of
$K_\Sigma^{1/2}\otimes L_\Sigma$.    Using the Riemann-Roch formula we deduce
\[
\dim H^0(K_\Sigma^{1/2}\otimes L_\Sigma) -\dim H^1(K_\Sigma^{1/2}\otimes L) = \deg L_\Sigma.
\]
On the other hand, Serre duality implies $\dim H^1(K_\Sigma^{1/2}\otimes
L)=H^0(K_\Sigma^{1/2}\otimes L_\Sigma^*)=h_{1/2}(L_\Sigma^*)$, where $L_\Sigma^*$ denotes
the dual of $L_\Sigma$. Hence
\be
h_{1/2}(L_\Sigma)-h_{1/2}(L^*)= \deg L_\Sigma.
\label{eq: sym}
\ee
Let $A$ and $B$ as in (\ref{eq: pull}). In [N] we proved the following result.

\begin{proposition} {\rm (a) For all $r\in(0,1]$ the splitting
${\bS}_L={\can}^{-1/2}\otimes L \, \oplus \, {\can}^{1/2}\otimes L$ induces a
splitting of $\ker D_{B,r}$ and moreover we have an isomorphism
\[
\ker D_{B,r}\cong H^0(K_\Sigma^{1/2}\otimes L_\Sigma^*) \oplus  H^0(K_\Sigma^{1/2}\otimes L_\Sigma)
\] 
so that $\dim \ker D_{B,r}=h_{1/2}(L_\Sigma^*)+h_{1/2}(L_\Sigma)$.

\noindent (b) There exist $r_0 >0$ and $z_0>0$ such that for all $r\in (0,r_0]$ the only
eigenvalue of $D_{B,r}$ in $(-z_0,z_0)$ is $0$.}
\label{prop: 210}
\end{proposition}

With respect to the splitting ${\bS}_L={\can}^{-1/2}\otimes L \, \oplus \, {\can}^{1/2}\otimes L$ 
the operator  $D_{B,r}$ has the block decomposition $D_{B,r}=Z_{B,r} + T_B$
where
\[
Z_{B,r}=\left[
\begin{array}{cc}
{\bf i}\nabla^{B}_{\zeta_r} & 0 \\
0  & -{\bf i}\nabla^{B}_{\zeta_r}
\end{array}
\right]
\]
and
\[
T_B=\left[
\begin{array}{cc}
0 & {}^\flat\bar{\partial}_B \\
{}^\flat\bar{\partial}_B^* & 0 
\end{array}
\right].
\]
Also ${\dir}_{B,r}=D_{B,r}+\lambda_r/2$.   Another important piece of
information is a supercommutator identity established in [N]. In our special
case it has the form
\be
\{Z_{B,r}, T_B\}:= Z_{B,r}T_B +T_B Z_{B,r} =0.
\label{eq: scom}
\ee

\noindent {\bf Step 1}  The same argument as in [Z]  proves the following
result.
\begin{proposition}
\[
\eta_0:=\lim_{r\ra 0}\eta({\dir}_{B,r})= -2\int_\Sigma \hat{\bf
A}(\nabla^\Sigma)\cdot {\bf ch}(B)\cdot
\frac{\tanh (c/2)-c/2}{c\tanh(c/2)} +\sum_\mu {\rm sign}\, \mu
\]
{\rm where ${\bf ch}(B)$ denotes the Chern-Character defined in terms of the
connection $B$ on $L\ra \Sigma$ and the remaining terms have the same
significance as in Theorem \ref{th: z}.}
\label{prop: z1}
\end{proposition}

Proceeding exactly as in \S 2.2  we conclude  (via Proposition \ref{prop: 210})
that
\[
\eta_0=\ell/6  -{\rm sign}\,(\ell) (\,h_{1/2}(L_\Sigma^*)+h_{1/2}(L_\Sigma)\,)
\]
 If we set $\xi_r=\xi({\dir}_{B,r})=\frac{1}{2}\eta_r$ we deduce
\[
\xi_0-\xi_r=\ell/12 -{\rm sign}\,(\ell)\cdot
\frac{h_{1/2}(L_\Sigma^*)+h_{1/2}(L_\Sigma)}{2}
-\xi_r.
\]
On the other hand, we have
\[
\xi_0-\xi_r=-\frac{\ell}{12}(\ell^2r^4-\chi r^4).
\]
This follows  from the first transgression formula.   We can quote this formula
since as $r\ra 0$ the only constituent of ${\dir}_{B,r}$ that changes is the
Levi-Civita connection while the coupling connection is independent of $r$.  The
 degree 3 part of the transgression of the index density $\hat{\bf A}
\wedge {\bf ch}(B)$    equals precisely the transgression of  the $\hat{\bf
A}$-genus which was computed in the first transgression formula. We conclude
\be
\frac{1}{2}\eta_r=\xi_r({\dir}_{B,r})=\frac{\ell}{12}-{\rm
sign}\,(\ell)\cdot \frac{h_{1/2}(L_\Sigma^*)+h_{1/2}(L_\Sigma)}{2}
+ \frac{\ell}{12}(\ell^2r^4-\chi r^2). 
\label{eq: eta1}
\ee
\noindent {\bf Step 2}   Now we ``transgress'' from ${\dir}_{B,r}$ to $D_{B,r}$
using the same deformation $(\tilde{\dir}_{B,r,\tau})$ as in \S 2.3. As in that
case we have
\be
\xi(D_{B,r})=\xi({\dir}_{B,r})+ SF(\tilde{\dir}_{B,r,\tau};0\leq \tau \leq
1) +\int_NT\hat{\bf A}(\nabla^{r,0}, \nabla^{r,1}).
\label{eq: deform}
\ee
Again there is no transgression term coming from the coupling which does not
change as $\tau$ runs from $0$ to $1$.

The spectral flow contribution occurs only at $\tau=1$ and using Proposition
\ref{prop: 210} we  determine it to be
\be
SF(\tilde{\dir}_{B,r,\tau};0\leq \tau \leq 1)=\left\{
\begin{array}{rll}
h_{1/2}(L_\Sigma^*)+h_{1/2}(L_\Sigma)&,& \ell>0 \\
0&,& \ell<0
\end{array}
\right.
\label{eq: spfl}
\ee
The transgression term is given by the second transgression formula. Putting
all the above together we  deduce
\be
\xi(D_{B,r})=\frac{\ell}{12}+\frac{h_{1/2}(L_\Sigma^*)+h_{1/2}(L_\Sigma)}{2}.
\label{eq: eta2}
\ee
Note that 
\be
\eta(D_{B,r})= \frac{\ell}{6}.
\label{eq: adiaeta}
\ee
Surprisingly, $\eta(D_{B,r})$ carries very little geometrical  information.    The extreme generality of the Bismut-Cheeger-Dai theorem may obscure   some beautiful symmetries  responsible for (\ref{eq: adiaeta}). We refer the reader to Appendix C where  we have included an elementary derivation of  this equality which works in the more general context of  Seifert manifolds and  we belive contains several illuminating informations.

\medskip 

\noindent{\bf Step 3}     Finally we compute $\xi(D_{A,r})$. {\em  In the remaining part we will assume} $k\in {\bZ}\cap
(0, |\ell|)$.  Note that if $k=0$ then $D_{A,r}=D_{B,r}$ and there is nothing to compute in this case.  Hence we have to consider only the case $0<k<\ell$.

The equality $A=B-\frac{k{\bf i}}{\ell} \vfi$ suggests using the path of connections
\[
B_t=B+t\dot{B},\;\;\dot{B}=-\frac{k{\bf i}}{\ell}\vfi_r.
\]
We have (omitting the $r$-subscript for brevity)
\be
\xi(D_A)=\xi(D_B) +SF(D_{B_t}) +\int_N T(\hat{\bf A} \wedge {\bf ch})(B_1,B_0).
\label{eq: var3}
\ee
This time $\hat{\bf A}$ is fixed and only the coupling connection changes. We have the
following result.
\begin{lemma}{\bf (Third transgression formula)}
\[
\int_N T(\hat{\bf A} \wedge {\bf ch})(B_1,B_0)=\frac{k^2}{2\ell}.
\]
\end{lemma}

\noindent {\bf Proof of the lemma} \hspace{.3cm} As we mentioned before,  the
only part which contributes to the transgression is ${\bf ch}$ through its
degree 4 component $\frac{c_1^2(B_t)}{2}$. For an arbitrary connection
$\nabla$ on $L$ we have
\[
\frac{c_1^2(\nabla)}{2}=-\frac{1}{8\pi^2} F(\nabla)\wedge F(\nabla).
\]
Thus in our case the transgression is
\[
T{\bf ch}= -\frac{1}{4\pi^2}\int_0^1 \dot{B} \wedge F_{B_t}.
\]
A simple computation shows
\[
F_{B_t}=F_{B} +td\dot{B} =F_B -\frac{tk {\bf i}}{\ell}d \vfi =F_B+2tk{\bf
i}\vfi^1\wedge \vfi^2
\]
\[
\dot{B} \wedge F_{B_t}=-\frac{k{\bf i}}{\ell}\vfi \wedge F_B
+\frac{2tk^2}{\ell} \vfi\wedge \vfi^1\wedge \vfi^2.
\]
Hence
\[
T{\bf ch}=\frac{{\bf i}k }{4\pi^2 \ell}\vfi \wedge F_B-\frac{k^2}{4\pi^2\ell} \vfi\wedge \vfi^1\wedge \vfi^2
\]
\[
 =\frac{k}{2\pi \ell}\vfi \wedge \frac{{\bf i}F_B}{2\pi} -\frac{k^2}{4\pi^2\ell}\vfi\wedge \vfi^1\wedge \vfi^2 .
\]
The lemma follows  integrating over $N$ and using the equalities
\[
\int_{base}\frac{{\bf i}F_B}{2\pi}=\deg L =k,\;\;\int_{fiber}\vfi =2\pi, \;\;\int_N \vfi \wedge  \vfi^1\wedge \vfi^2 =2\pi^2.\;\;\Box
\]

\medskip

To compute  the spectral flow in  (\ref{eq: var3}) we need to  go deeper inside
the structure of $D_{B_t}$. We have
\[
D_{B_t}=Z_{B_t,r}+T_{B_t}
\]
Note that $T_{B_t}=T_B$ since $T$ involves only derivatives along horizontal
directions  while $B_t$ changes only in the vertical direction. As for $Z_{B_t,r}$
we have
\[
Z_{B_t,r}=Z_{B,r}=\left[
\begin{array}{cc}
{\bf i}\nabla^{B}_{\zeta_r}+{\bf i}t\dot{B}(\zeta_r) & 0 \\
0  & -{\bf i}\nabla^{B}_{\zeta_r}-{\bf i}t\dot{B}(\zeta_r)
\end{array}
\right]
\]
\[
=Z_{B,r}+\frac{t}{r}\left[
\begin{array}{cc} k/\ell & 0 \\
0 & -k/\ell 
\end{array}
\right]
\]
Denote the  ``matrix'' above by ${\gR}$  and set
$Z_{B}:=Z_{B,r=1}$. Using the equality $\zeta_r=r^{-1}\zeta$ we deduce
\[
Z_{B_t, r}=\frac{1}{r}(Z_B +t{\gR}).
\]
Observe now that both $Z_B$ and ${\gR}$ anticommute with $T_B$ so that
\[
D^2_{B_t,r}=\frac{1}{r^2}(Z_B+t{\gR})^2 +T_B^2.
\]
In particular, this shows
\[
\ker D_{B_t,r}=\ker D^2_{B_t,r}=\ker (Z_{B}+r{\gR})\, \cap \, \ker T_B.
\]
{\em Since} $ 0 < k < |\ell|$  we see that $\ker Z_B +t{\gR}=(0)$ if $t\in (0,1]$. In
other words, the only contribution to the spectral flow arises at $t=0$. 
Denote by $\{\mu_i(t)\}$ the eigenvalues of $D_{B_t,r}$ such that $\mu_i(0)=0$.
There are $\dim \ker D_{B,r}$ such eigenvalues and denote by $\sigma_-$ the
number of those such that $\dot{\mu}_i(0) <0$. The spectral flow is then
$-\sigma_-$.  Determining the eigenvalues $\mu_i(t)$ may be a complicated job.
We  follow a different description of $\sigma_-$ given in [RS].  

Set $E_0= \ker D_{B,r}$, denote by $P_0$ the orthogonal projection onto $E_0$
and define the {\em resonance matrix} $R:E_0\ra E_0$ by 
\[
R=P_0\dot{D}_{B_t,r}\!\mid_{t=0}:E_0\ra E_0.
\]
Clearly $R$ is nondegenerate and, as explained in [RS],  $\sigma_-$ can be identified with $\sigma_-(R)$
which is the number of negative eigenvalues of $R$.  This number  can be
determined using the explicit description of ${\gR}$ and Proposition \ref{prop:
210} (a). We deduce
\[
SF(D_{B_t,r})=-\sigma_- =\left\{
\begin{array}{rll}
-h_{1/2}(L_\Sigma)&,& \ell>0 \\
-h_{1/2}(L_\Sigma^*)&, & \ell <0
\end{array}
\right.
\]
Using the  third transgression formula and the  equalities (\ref{eq: var3}),  (\ref{eq: eta2}) we 
finally determine
\be
\xi(D_A)=-\frac{\ell}{24} +\frac{k^2}{2\ell} + {\rm
sign}\,(\ell)\cdot\frac{ h_{1/2}(L_\Sigma^*)-h_{1/2}(L_\Sigma)}{2} \stackrel{(\ref{eq:
sym})}{=}\frac{\ell}{12}+\frac{k^2}{2\ell}-{\rm sign}\,(\ell) \frac{k}{2}
\label{eq: eta3}
\ee
Again $\xi(D_{A,r})$ is a topological quantity!!!

\begin{remark}{\rm  The spectral flow computation in Step 3 used in an essential way the fact that $k\in (0,\ell)$. In fact, if  we started with a different $k'\equiv k $mod $\ell$ then the computations at  Step3 would be affected in  both the transgression term and in the spectral flow term ( which would now have several contributions).  One can verify easily that these changes  cancel each other so that the final result  is independent of the  choice of  a residue of $k$ mod $\ell$.}
\end{remark}

\section{Finite energy Seiberg-Witten monopoles}
\setcounter{equation}{0}

Throughout this section, a hat over an object  will 
signal (unless otherwise indicated) that  it  is a 4-dimensional geometric
object.  

For example, if $N$ is a 3-manifold then on the tube ${\bR}\times N$ there
exist two exterior  derivatives: the 3-dimensional exterior derivative $d$ along the
slices $\{t\}\times N$ and the 4-dimensional exterior derivative $\hat{d}$ so
that $\hat{d}=dt \wedge \partial_t + d$. If $A(t)$ is a family of connection on
some vector bundle $E\ra N$ then we get a bundle $\hat{E}\ra {\bR}\times N$ 
and  we can think of the path $A(t)$ as a  connection $\hA$ on $\hat{E}$.  We
will denote by $F_{A(t)}$ the curvature of $A(t)$ on the slice $\{t\}\times N$
while $\hat{F}_\hA$ will denote the curvature of $\hA$ on the tube.

\subsection{The 4-dimensional Seiberg-Witten equations}
 Let $\hat{N}$ denote an oriented  4-manifold ({\em not necessarily compact}), equipped with a Riemann metric $\hat{g}$. Denote by $\hat{\ast}$ the Hodge star operator
 defined by the metric $\hg$ and the orientation of $\hat{N}$.  {\em Fix} a connection $\nah$ on $T\hat{N}$ compatible with $\hg$. $\nah$ need not be  the Levi-Civita connection.

  Denote by  $Spin_c(\hat{N})$ the  collection of isomorphism classes of $spin^c$ structures on $\hat{N}$.  For each
  $\hsi \in Spin^c(\hat{N})$ we denote by  $\det \hsi$ the associated line bundle and by
 ${\bS}_\hsi ={\hbS}_\hsi^+\oplus {\hbS}_\hsi^-$ the associated bundle of spinors. 
 Note that $\det \hsi \cong \det \hat {\bS}_\hsi^+$.

 Denote by ${\gA}_\si$ the space of hermitian connections on ${\bS}_\si$
 compatible  with both the ${\bZ}_2$-grading and the fixed background
 connection $\nah$. More precisely, $A\in {\gA}_\hsi(\hat{N})$ if for any $\alpha \in
 \Omega^1(N)$, any $X\in {\rm Vect}\,(N)$ and any $\hpsi\in C^\infty({\bS}_\si)$ we have
 \[
 \nabla_X^A({\hbc}(\alpha)\hpsi )=\hbc(\nabla_X\alpha)\hpsi
 +\hbc(\alpha)\nabla^A_X\hpsi
 \]
 where 
 \[
 {\hbc}:T^*\hat{N}\ra {\rm Hom}\,({\hbS}_\hsi^+,
 {\bS}_\hsi^-)
 \]
 denotes the Clifford multiplication.  Any connection on $\det \hsi$ determines
 a connection in ${\gA}_\si$ and moreover, once we fix a connection $A_0\in
 {\gA}_\hsi(\hat{N})$,  we can identify ${\gA}_\hsi(\hat{N})$ with
 $\ii\Omega^1(\hat{N})$. To any
 connection $\hA\in {\gA}_\hsi(\hat{N})$ we can associate  the Dirac operator
 \[
 \hat{D}_\hA: \Gamma({\hbS}_\hsi^+)\ra \Gamma({\hbS}_\hsi^-)
 \]
 defined as the composition
 \[
 \Gamma({\hbS}_\hsi) \stackrel{\nah^\hA}{\ra} \Gamma(T^*\hat{N}\otimes
 {\hbS}_\hsi^+)\stackrel{\hbc}{\ra} \Gamma({\hbS}_\hsi^-).
 \]

 There is a natural quadratic map 
 \[
 q: \Gamma({\hbS}_\hsi^+)\ra {\rm End}\, ({\hbS}_\hsi^+),\;\;\hpsi \mapsto
 \tau(\hpsi)
 \]
 defined by
 \[
 q(\hpsi)\hat{\phi}=\lan \hat{\phi}, \hpsi\ran -\frac{1}{2} |\hpsi|^2
 \hat{\phi}.
 \]
 In terms of Dirac's bra-ket notation $\tau(\hpsi)$ can be alternatively
 described as 
 \[
 q(\lan\hpsi|)=|\hpsi\ran \lan \hpsi | -\frac{1}{2}\lan \hpsi
 | \hpsi\ran.
 \]
 Note  that for each $\hpsi$ the endomorphism $\tau(\hpsi)$ is symmetric and
 traceless (see Appendix $D$). 
 
 The quantization map from the exterior algebra to the Clifford algebra extends
 the Clifford multiplication to a map 
 \[
 \hbc : \Lambda^* T^*\hat{N} \ra  {\rm End}\, ({\hbS}_\hsi)
 \]
This map has the property that $\hbc(\omega)$ is a traceless, skew-symmetric
endomorphism of ${\hbS}_\hsi^+$ for any  $\hg$-self-dual real valued 2-form
$\omega$.

The Seiberg-Witten equations  (associated to the $spin^c$ structure $\hsi$) are
equations for a pair $(\hpsi, \hA)$ = (spinor in ${\bS}_\hsi^+$, connection in
${\gA}_\hsi(\hat{N})$).  More precisely, they are
\[
(\widehat{SW})\;\;\left\{ 
\begin{array}{rcl}
\hat{D}_\hA \hpsi & = & 0 \\
\hbc(\hat{F}^+_\hA)& = & \tau(\hpsi) 
\end{array}
\right.
\]

In the remaining part of this subsection we will make further  additional
assumptions on the geometry and the topology of $\hat{N}$ and explain how this
affects the Seiberg-Witten equations.

\begin{figure}
\centerline{\psfig{figure=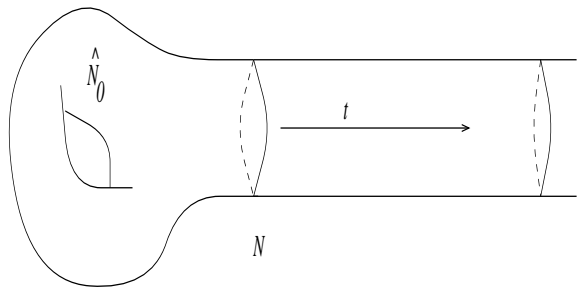,height=1.5in,width=3in}}
\caption{\sl The background manifold $\hat{N}$}
\label{fig: cyl1}
\end{figure}

More precisely, assume the manifold $\hat{N}$ can be decomposed as
\[
\hat{N}=\hat{N}_0 \cup [0,\infty) \times N
\]
where $\hat{N}_0$ is a compact oriented 4-manifold with boundary $\partial
\hat{N}_0 =N$.   We will denote  by $t$ the longitudinal coordinate on the
cylindrical part of $N$ (see Figure \ref{fig: cyl1}).
Fix a tubular neighborhood  $(-1,0]\times N$ of $N$ in $\hat{N}_0$, a metric
$g$ on $N$ and a connection $\nabla$ compatible with $g$, {\em not necessarily the Levi-Civita
connection} of $g$. We assume that along the  infinite cylinder $(-1,\infty)\times N$
 the metric $\hg$ is a  product metric $\hg =dt^2 + g$. We fix a   connection
 $\nah$ compatible with $\hg$ such that along the above cylindrical end  it has
 the form
 \[
 \nah = \pat \wedge dt + \na .
 \]
 We denoted by $\pat$ the  $\hg$-gradient of $\tau$ where $\tau :\hat{N} \ra [0,\infty)$ is a smooth function
  which coincides with the  canonical projection $[0,\infty)\times N \ra
  [0,\infty)$ on the infinite neck.

  Note that the  $spin^c$ structure $\hsi$ induces  a  $spin^c$
  structure $\si$ on  $N=\partial \hat{N}_0$.  Denote by ${\bS}_\si\ra N$ the
  associated  bundle of spinors and by $\bc :T^*N\ra {\rm End}\, (N)$ the
  corresponding Clifford multiplication. As in the $4$-dimensional case we can
  define ${\gA}_\si(N)$.
  
  Fix a reference connection $\hA_0\in {\gA}_\hsi(\hat{N})$ which  along the neck is
  gauge equivalent to a product connection $dt\otimes \partial_t +A_0$, $A_0\in
  {\gA}_\si(N)$.   Now define the configuration space $\hconf$ as the set of pairs 
  $(\psi,  \hA_0+\ii \hat{a}):=(\hpsi, \hA)$=(spinor, connection) such that    
  $(\hpsi, \ii \hat{a}) \in L^{2,2}_{loc}(\hat{\bS}\oplus{\bf  i}T^*\hat{N})$ and  
  \[
\nah^{\hA}_\pat\hpsi \oplus i_\pat\left(\ii\hat{F}_{\hA}\right) \in
  L^2(\hat{\bS}_\hsi \oplus  \ii\Lambda^1T^*\hat{N}).
  \] 
 We denoted    by $i_\pat$  the contraction by $\pat$.  For brevity,  will denote the elements of $\hconf$ by   the generic symbol $\hco$. 

  \begin{definition}{\rm  (a) A finite energy solution of $(\widehat{SW}_\omega)$
  is a solution $(\hpsi, \hA)$ such that $(\hpsi, \hA-\hA_0)\in \hconf$. 
  
  \noindent (b) A Seiberg-Witten tunneling is a finite energy solution on
  $\hat{N}={\bR}\times N$.}
  \end{definition}
  
  There is an infinite dimensional group $\hgauge$ acting on the configuration space, more precisely 
  \[
  \hgauge =\{ \gamma \in  {\rm Map}\,(\hat{N}, S^1)\; ; \;  \gamma \in
  L^{3,2}_{loc}\}  
  \]
 The group $\hgauge$ acts (on the right) on $\hconf$ and  transforms finite
  energy solutions  to finite energy solutions. Define
  \[
  \hmodu :=\{ (\hpsi, \hA) \; {\rm finite\; energy\; solutions\; of\;}
  \widehat{SW}\}/\hgauge.
  \]
  In this section we want to analyze the the Fredholm properties of the deformation complex naturally associated to 
  $\hmodu$ when $N$ is a circle bundle over a Riemann surface.  In particular
  we will compute the virtual dimension of the space of Seiberg-Witten
  tunnelings. 
  
  We conclude this subsection  with  a simple but crucial observation which
  reveals the dynamical feature of the  Seiberg-Witten equations on cylinders
  which perhaps will explain the tunneling terminology.

   Note that if we set $J=\hbc(d\tau)$
  then  $J$ induces  isomorphisms
  \be
  \hat{\bS}^+_\hsi\!\mid_N\cong \hat{\bS}_\hsi^-\!\mid_N \cong {\bS}_\si
  \label{eq: spin1}
  \ee
  and
  \be
  \bc(\alpha)=J\hbc(\alpha),\;\; \forall\alpha \in \Omega^1(N) \hookrightarrow \Omega^1 ([0,\infty)\times
  N).
  \label{eq: spin2}
  \ee
 A connection $\hA\in {\gA}_\hsi(\hat{N})$ is said to be in a {\em temporal gauge} if
 $i_\pat (\hA- \hA_0)=0$ along  the infinite neck $[0,\infty) \times N$.

 Assume now that $(\hpsi, \hA)$   is a finite energy solution of
 $(\widehat{SW})$ such that $\hA$ is in a temporal gauge.  Along the neck we
 can write
 \[
 \hpsi =\psi(t),\;\; \hA = A_0 + \ii a(t)
 \]
 where $A_0=\hA_0\!\mid_N$, $\psi(t) \in \Gamma ({\bS}_\si)$, $a(t) \in  \Omega^1(N)$, $\forall t\geq 0$.  Then  (along the neck)
 \be
 \hat{F}^+_{\hA} =\frac{1}{2}\{ (F_a +\ast \ii \dot{a}) +dt \wedge (\ii\dot{a} +\ast
 F_a)\}
 \label{eq: sd}
 \ee
 where $A_0 + a(t)$ is the connection on  the line bundle $\det \si$ restricted to  the slice
 $\{t\}\times N$,  $F_a =F_{A_0 +\ii a}$ denotes is curvature and $\ast$ denotes the Hodge star operator on $N$. $A_0+\ii a(t)$  induces a  Dirac operator 
 \[
 {D}_a={D}_{a(t)}: \Gamma({\bS}_\si)\ra \Gamma({\bS}_\si).
 \]
 Using  (\ref{eq: spin1}) and (\ref{eq: spin2}) we deduce that along the neck
 \[
 \hat{D}_\hA= J\left(\partial_t -{\dir}_a\right).
 \]
 The equality (\ref{eq: sd})   now implies
 \[
 \hbc(\hat{F}^+_\hA)= \bc (\ast F_a + \ii \dot{a}).
 \]
 Consequently, along the neck, in a temporal gauge,  the Seiberg-Witten equations
 can be rewritten as
 \be
 \left\{\begin{array}{rcl}
 \dot{\psi} &=& {D}_a \psi \\
 \ii\bc (\dot{a})&=&q(\psi)-\bc(\ast F_a) 
 \end{array}
 \right. .
 \label{eq: swtube}
 \ee
 The right-hand-side  of (\ref{eq: swtube})    arises  when one considers the three dimensional counterpart of the Seiberg-Witten equations.

\subsection{The 3-dimensional Seiberg-Witten equations}
To formulate these equations we need to consider a new configuration space. Fix
a connection $A_0\in {\gA}_\si(N)$ and define
\[
\conf =\{(\psi, A)\; ;\; (\psi, (A-A_0) \in L^{1,2}({\bS}_\si \oplus {\bf i}T^*N)\}.
\]
  For brevity, its elements will be denoted by  the symbol $\co$ and we will often write
  $\co=(\psi, a)$ instead of $(\psi, A_0+\ii a)$ whenever no confusion is possible. There is an energy functional  ${\en} :\conf \ra {\bR}$ defined by
\be
{\en}(\psi, A)=\frac{\ii }{2}\int_N a\wedge (F_{A_0}+F_A) +\frac{1}{2}\int_N\lan \psi ,{\dir}_A\psi\ran dv_g.
\label{eq: energy}
\ee
The gauge group 
\[
\gauge =\{\gamma\in {\rm Map}\,(N, S^1)\; ;\; \gamma \in L^{2,2}\}
\]
acts on $\conf$ and moreover
\[
{\en}(\gamma^{-1}\cdot(\psi,A))-{\en}(\psi, A)=-\int_N\gamma^{-1}d\gamma \wedge F_{A_0} =2\pi{\bf
i}\int_N\gamma^{-1}d\gamma \wedge c_1(A_0)
\]
where  we denoted by $c_1(A_0)$ the  2-from representing the first  Chern class of $\det \si$  associated to $A_0$ via the Chern-Weil construction.  The  $L^2$-gradient  of ${\en}$  is (see
[N] or [MOY])
\[
\nabla {\en}(\psi, A) =\left[
\begin{array}{c}
{\dir}_A\psi \\
q(\psi)-\ast F_A
\end{array}
\right]
\]
where  we  tacitly  identified $q(\psi)$ with a purely imaginary  1-form via the Clifford multiplication.   The 3-dimensional Seiberg-Witten equations  can now be  described  as
\[
\nabla{\en} (\co)=0 \Longleftrightarrow \left\{
\begin{array}{rcl}
{\dir}_A\psi & =& 0 \\
{\bc}(\ast F_A) & =& q(\psi)
\end{array}
\right.
\]
We see that  (\ref{eq: swtube}) can be rewritten as a gradient flow equation 
\be
\dot{\co}=\nabla {\en}(\co).
\label{eq: swtube1}
\ee
This last equation suggests that  as $t\ra \infty$  $\co(t)$ converges to a critical point  of ${\en}$.   Assuming the finite energy  condition this  can be proved for arbitrary $N$ using the techniques of [MMR].  However, unlike  the Yang-Mills situation,  the nature of critical points  and  the manner in which they are organized are less transparent in the Seiberg-Witten case.  This is  the reason
 why   we will concentrate on a  special case.

\medskip

{\em In remaining  of the section $N$ will be assumed to be  a degree $\ell \neq 0$  circle bundle over a  genus $g>0$ Riemann surface $\Sigma$,  $S^1\hookrightarrow N \stackrel{\pi}{\ra} \Sigma$ equipped with  the metric described in $\S 2.1$  As background $g$-compatible connection on $N$ we choose the adiabatic connection $\nabla^0=\lim_{t\ra 0}\nabla^{r,t}$, where $r$ is fixed and small.}

\medskip

The $spin^c$ structures on $N$ are bijectively parameterized by   the space
of isomorphism classes of hermitian line bundles on $N$.  Fix a $spin$
structure on $\Sigma$ determined by a holomorphic square root $K^{1/2}$. If  $L\ra N$ is such a line bundle then the corresponding bundle of complex spinors is
\[
{\bS}_L={\can}^{-1/2}\otimes L \oplus {\can}^{1/2}\otimes L.
\] 
 Moreover we can identify the connections in ${\gA}_\si(N)$ with the hermitian
 connections on $L$. The Dirac operator  on ${\bS}_L$  induced by $\nabla^0$ and a connection $A$
 on $L$  will be  denoted by $D_A$. If instead of $\nabla^0$ we use the 
 Levi-Civita connection of the metric $g_r$ we get a different Dirac type
 operator we denote by ${\dir}_A$.  The  operator $D_A$  can
 be related to ${\dir}_A$ by  the following simple identity (see Sect. \S 2.1).
 \[
 D_A= {\dir}_A -\lambda_r/2,\;\;\lambda_r =-r\ell.
 \]
 Both $D_A$ and $\dir_A$ have obvious extensions to $[0,\infty)\times N$  given
 by
 \[
 \hat{D}_A=J(\partial_t -D_A),\;\;\; \hat{\dir}_A=J(\partial_t -\dir_A).
 \]
 Under these special geometric circumstances the Seiberg-Witten equations can be rewritten in a more
 useful form.

 Using the decomposition ${\bS}_L\cong \left({\can}^{-1/2}\otimes L\right)\, \oplus
 \, \left({\can}^{1/2}\otimes L\right)$ we can  represent any section $\psi$ of ${\bS}_L$  as
 $\psi=\psi_- \oplus \psi_+$. Then the Seiberg-Witten equations can be  rephrased as  (see [N])
\be
\left\{
\begin{array}{rrrcl}
{\bf i}{\nabla}^A_\zeta \psi_- & +  \bar{\partial}_A \psi_+ & +\lambda  \psi_- &  = & 0 \\
 & & & & \\
(\bar{\partial}_A)^* \psi_- &  - {\bf i}{\nabla}_\zeta \psi_+  & +\lambda  \psi_+ & = &0 \\
&&&&\\

& &  \frac{1}{2}(|\psi_-|^2-|\psi_+|^2) & = &{\bf i}F_A(\zeta_1, \zeta_2) \\
& & & & \\
& & {\bf i}\psi_- \bar{\psi}_+&= &\bar{\ve}\otimes F_A(\zeta_1
+{\bf i} \zeta_2, \zeta)
\end{array}
\right.
\label{eq: sw}
\ee
where ${\ve} =2^{-1/2}(\vfi^1+\ii\vfi^2)$.

Set
\[
\conf^*=\{(\psi, A)\in \conf\; ;\; \psi \not\equiv 0\}.
\]
The configurations in  $\conf^*$ are called {\em irreducible}. As in [M]  one  
can show that ${\gB}:=\conf/\gauge$ is a metric space and, moreover,  ${\gB}^* =\conf^*/\gauge$ is  a Banach manifold.  
This is proved  using the  existence of local slices for the $\gauge$-action exactly as in the Yang-Mills case.    For every configuration $\co\in {\conf}$ we will denote by $[\co]$ its image in ${\gB}$.

The solutions of (\ref{eq: sw}) are explicitly described in [N] and [MOY]. 
Here are the relevant facts.

\medskip

\noindent {\bf Fact 1.}  If $c_1(L)$ is not torsion then  (\ref{eq: sw}) has no
solutions.

Assume now that $c_1(L) \equiv \kappa$ (mod $\ell$) and define
\[
R_\kappa=\{n\in {\bZ} \; ; \; 1\leq |n|\leq g-1, \;\; n\equiv \kappa \; {\rm mod}\; \ell\}.
\]

\noindent{\bf Fact 2.}  Any irreducible solution $(\psi, A)$   of (\ref{eq:
sw})  is gauge equivalent  to the pullback of a pair $(\tilde{\psi}, B)$
where $B$ is a connection in a line bundle $L_\Sigma \ra \Sigma$ such that
$\deg L_\Sigma \in R_\kappa$ (so that $\pi^*L_\Sigma \cong L$),
$\tilde{\psi}=\tilde{\psi}_-\oplus \tilde{\psi}_-$ is a
section of $C^\infty(K^{-1/2}\otimes L_\Sigma \, \oplus \, K^{1/2}\otimes L)$.
The connection $B$ defines holomorphic structures in $K^{\pm 1/2}\otimes L$.
$\tilde{\psi}_-$ is an antiholomorphic section of $K^{-1/2}\otimes L$  while
$\tilde{\psi}_+$ is a holomorphic section of $K^{1/2}\otimes L_\Sigma$.    Moreover,  one of
$\tilde{\psi}_-$ or $\tilde{\psi}_+$ is zero and  satisfy the identity:
\[
\frac{1}{4\pi}\int_\Sigma(|\tilde{\psi}_-|^2 -|\tilde{\psi}_+|^2) \, dv =\deg L_\Sigma.
\]
Thus $\tilde{\psi}_+ =0$ if $\deg L_\Sigma >0$ and $\tilde{\psi}_- =0$ if  $\deg L_\Sigma <0$.
The irreducible  part (mod $\gauge$), denoted by ${\modu}^*$  consists of $\# R_\kappa$ components
\[
{\modu}^*=\bigcup_{n\in R_\kappa} {\modu}_{\kappa, n}.
\]
  The component ${\modu}_n={\modu}_{\kappa,n}$ (corresponds to the choice $\deg L_\Sigma =n$) is diffeomorphic to a symmetric product
of $(g-1)-|n|$ copies of $\Sigma$ and thus has real dimension $2(g-1-|n|)$.
Each component is Bott nondegenerate as a critical set.
(Pairs $(\tilde{\psi}_-\oplus \tilde{\psi}_+, B)$ as above are known as {\em vortex
pairs} on $\Sigma$.)

\medskip

\noindent {\bf Fact 3.}  The reducible solutions  consist of pairs (zero
spinor, flat connection).   Modulo $\gauge$ they form a space ${\modu}_\kappa^0$ homeomorphic to  a $2g$-dimensional torus.  Moreover if $\kappa \not\equiv 0$ (mod $\ell$) the reducible part is nondegenerate (in a sense described in [MOY]). If $\kappa \equiv 0$  the reducible solutions  can be identified with the thetadivisor $W_{g-1}$ inside the jacobian $J_{g-1}(\Sigma)$ (see [GH] for a  definition of $W_{g-1}$).

\bigskip

Associated to each component ${\modu}$ there is a deformation theory
which we now proceed to describe.  We will concentrate only on the irreducible
part $\conf^*$. Since $c_1(L)$ is torsion the energy
functional ${\en}$ is gauge invariant  and  thus it descends to a well defined
functional 
\[
{\uen}: {\gB}^*\ra {\bR}.
\]
The group $\gauge$ is a Hilbert-Lie group and its Lie algebra can be identified
with the space ${\gog}:=L^{2,2}(N, \ii{\bR})$. The exponential map has the form
\[
{\gog}\ni \ii f \mapsto (\exp(\ii f): N \ra S^1).
\]
The tangent space to the orbit ${\cal O}_{\phi, A}$ through  $\co=(\phi, A)$ of the {\em right} action of
$\gauge$  is the range of the infinitesimal action
operator
\[
{\Lie}={\Lie}_{\co}:{\gog}\ra {\X}:=L^{1,2}({\bS}_L)\oplus L^{1,2}(\ii
T^*N),\;\;\ii f\mapsto -\ii f \oplus \ii df.
\]
The tangent space to ${\gB}^*$ at  $[\co]$ can be identified with the orthogonal complement to the tangent space to the orbit ${\cal O}_\co$ and ultimately with the kernel of ${\Lie}^*_\co$, the adjoint of  ${\Lie}_\co$. An integration by parts shows
\[
{\Lie}^*(\dps\oplus \ida)=-\ii d^* a+\ii \im \lan \phi, \dps\ran, \;\;\forall \dps\oplus \ida \in {\X}.
\]
We can use the affine structure of ${\conf}$ to linearize  $\nabla {\en}$ at a given configuration $\co=(\phi, A)$ and we  obtain the {\em unrestricted hessian} at $\co$
\[
{\hes}_\co\left[
\begin{array}{c}
\dps \\
\ida
\end{array}
\right] =\frac{d}{dt}\!\mid_{t=0}\nabla{\en}(\psi +t\dps, A+t\ida)= \left[
\begin{array}{c}
D_A{\dps} +{\bf c}({\ida})\phi  \\
-{\bf i}\ast d{\da} +{\dta}(\phi, {\dps}) 
\end{array}
\right]
\]
The term $ {\dta}(\phi, {\dps}) $   is formally defined by the equality
\[
{\dta}(\phi, {\dps})  :={\dt}\!\mid_{t=0}q(\phi +t{\dps})
\]
where we regard $q$ as a quadratic map $q:{\bS}_L \ra {\bf i}T^*N$.   

The {\em stabilized hessian} of ${\en}$ at $\co=(\phi, A)$ is the unbounded operator on $L^2({\bS}_L\oplus\ii(\Lambda^1\oplus\Lambda^0)T^*N)$ defined by
\[
{\hhes}_{\co}\left[
\begin{array}{c}
{\dps}\oplus{\ida} \\
{\bf i}f 
\end{array}
\right ]:=\left[\begin {array}{cc}
{\hes} & {\Lie} \\
{\Lie}^* & 0
\end{array}
\right]  \left[
\begin{array}{c}
{\dps}\oplus {\ida} \\
{\bf i}f 
\end{array}
\right ] = \left[
\begin{array}{clcl}
D_A\dps  &                               &+ &  {\bc}(\ida)\phi-\ii f\phi   \\
                 &-\ii\ast  d\da +\ii df & +&\dta(\phi, \dps) \\
                  &\ii d^*\da& + & \ii \im \lan \phi, \dps\ran  
\end{array}
\right]
\]
In  [N] and [MOY]  it is shown that if $[\co]\in {\modu}_{\kappa, n}$ then  the kernel of the stabilized hessian ${\hhes}_\co$  is naturally isomorphic to the tangent space $T_{[\co]}{\modu}_{\kappa,n}$.  Now define 
\[
{\hhes}_0  \left[
\begin{array}{c}
{\dps}
\\
 {\ida} \\
{\bf i}f 
\end{array}
\right ] =\left[
\begin{array}{crc}
D_A & 0 & 0\\
0 & -\ast d & d \\
0 & d^*  & 0 
\end{array}
\right]  \left[
\begin{array}{c}
{\dps}
\\
 {\ida} \\
{\bf i}f 
\end{array}
\right ] 
\]
and ${\cal P}={\cal P}_\phi$ by
\[
{\cal P}_\phi \left[
\begin{array}{c}
{\dps}
\\
 {\ida} \\
{\bf i}f 
\end{array}
\right ] 
=\left[
\begin{array}{c}
{\bc}(\ida)\phi -\ii f \phi \\
\dta(\phi, \dps) \\
\ii \im \lan \phi, \dps\ran
\end{array}
\right]
\]
Note that  ${\hhes}_0={\hhes}_0(\co)$ is an elliptic selfadjoint operator  for any $\co \in {\conf}$ and ${\hhes}_\co={\hhes}_0+{\cal
P}_\phi$.     For every $\co \in {\conf}$ define $SF_\pm(\co)\in {\bZ}$ as the spectral flow of the path $\pm({\hhes}_0(\co) +t{\cal
P}_\phi)$, $t\in [0,1]$. The next subsection is devoted to the computation of $SF(\co)$ when $[\co]\in {\modu}_{\kappa, n}$. 
For the reducible component ${\cal P}_\phi\equiv 0$ and this problem is
trivial.

Define now  for later use the {\em resonance matrix}. This is the quadratic form ${\cal R}$ on $\ker {\hhes}_0$ defined by
\[
{\cal R}\Xi={\cal R}_\phi\Xi={\bf Proj}{\cal P}_\phi\Xi ,\;\;\; \Xi= \dps\oplus \ida \oplus \ii f \in \ker {\hhes}_0
\]
where {\bf Proj} denotes the orthogonal projection onto $\ker {\hhes}_0$. Note  also that  for $[\co]\in {\modu}_{\kappa, n}$
\[
\ker {\hhes}_0([\co]) = H^0(K^{1/2}-L_\Sigma)\oplus H^0(K^{1/2}+L_\Sigma) \oplus H^1(\Sigma, {\bR})\oplus H^0(\Sigma, {\bR})
\]
where $L_\Sigma$ is the holomorphic line bundle on $\Sigma$ determined by  $[\co]$ as in  {\bf Fact 2}.

\subsection{Spectral flows and perturbation theory}
Fix $[\co]=[\phi, A]\in {\modu}_{\kappa, n}$. Assume for simplicity that $n<0$ so that $\phi_-=0$. Denote by $L_\Sigma$ holomorphic line bundle on $\Sigma$
($\deg L_\Sigma =n$), $B$ the induced connection on $L_\Sigma$   and by $\tilde{\phi}_+$ the holomorphic section of $K^{1/2}\otimes L_\Sigma$ determined by $[\co]$. The computation of $SF_+([\co])$ is carried  in two steps.  We consider only the  spectral flow $SF_+$. Also for simplicity we will  write $L$ instead of
$L_\Sigma$, and $\phi$ instead of $\tilde{\phi}$.

\medskip

\noindent{\bf Step 1} Along the path $t\mapsto {\hhes}_t:={\hhes}_0+t{\cal P}$,  $t\in [0,1]$ there is no spectral flow contribution for $t\neq 0$.

\medskip

\noindent  {\bf Step 2}   Compute the spectral flow  contribution at $t=0$.

\bigskip

Note first that $t\mapsto {\hhes}_t$ is an analytic family of selfadjoint operators with compact resolvent and thus  by known perturbation results  (see
[Kato], Thm. 3.9, Chap. VII) the eigenvalues and the eigenvectors of this family can be locally organized in  analytic families.  To complete  the first  step it suffices to show that $\dim \ker{\hhes}_t$  is independent of $0<t\leq 1$.   

To this aim consider  as in [N], Sec. 4.2, the following elliptic complex
\[
(V_{[\co]}):\;\;0 \ra{\bf i}\Omega^0(\Sigma)\stackrel{I}{\ra}\Gamma(L\otimes K^{1/2})\oplus {\bf i}\Omega^1(\Sigma) \stackrel{\Upsilon}{\ra}\Gamma(L\otimes K^{-1/2})\oplus {\bf i}\Omega^0(\Sigma) \ra 0.
\]
where
\[
\Upsilon
\left[
\begin{array}{c}
\dot{\beta}\\
{\ida}
\end{array}
\right]= \left[
\begin{array}{c}
\bar{\partial}\dot{\beta} +2^{-1/2}{\ida}^{0,1}\tilde{\phi}_+ \\
{\bf i}\ast d{\da}-{\bf i}{\re}\, \lan \tilde{\phi}_+, \dot{\beta}\ran
\end{array}
\right]
\]
${\ida}^{0,1}$ component is the $K^{-1}$-component of ${\ida}$ corresponding to the  orthogonal decomposition $T^*\Sigma \otimes {\bC}\cong K\oplus K^{-1}$. $I$ is the infinitesimal action 
\[
\ii f \stackrel{I}{\mapsto} (-{\bf i} f \tilde{\phi}_+ , B+{\bf i}df).
\]
In Sec. 4.2 of  [N] it is shown that  $H^0(V_{[\co]})\cong H^2(V_{[{\co}]})\cong 0$ and
\[
\dim_{\bf R} H^1(V_{[\co]})= -{\rm ind}_{\bf R}\,  (V_{[\co]})=\dim {\modu}_{\kappa, n}.
\]
Arguing exactly as in Sec. 5.6 of [MOY] one can prove that 
\[
\ker {\hhes}_t\cong H^1(V_{[\co]}),\;\;\forall t\in (0,1].
\]
In particular, if $\Xi=\dps \oplus \ii\dot{a} \oplus \ii f \in \ker {\hhes}_t$,
$t>0$ then $f\equiv 0$. This concludes  the first step in our program.

\bigskip

{\bf Step 2} Before we embark in the  computation of the spectral flow  contribution at $t=0$ we
need to  survey a few facts pertaining to  perturbation theory.

As we have already mentioned, the spectral data of ${\hhes}_t$ can be organized
in families  depending analytically upon $t$. Denote by $Z$ the set of all pairs
$(\lambda(t), \Xi(t))$ where $\lambda(t)$ is an eigenvalue of ${\hhes}_t$,
$\Xi(t)$  is a (length 1) eigenvector  corresponding to $\lambda(t)$,
$\lambda(0)= 0$ and the dependence
\[
t\mapsto (\lambda(t), \Xi(t))
\]
is analytic. Clearly, $\#Z=\dim {\hhes}_0$. For every $(\lambda(t), \Xi(t))$ we
have Taylor expansions
\[
\lambda(t)=\lambda_\nu t^\nu+ \cdots ,\;\;\lambda_\nu \neq 0
\]
\[
\Xi(t)=\Xi_0+t\Xi_1+\cdots ,\;\;\Xi_0\in \ker {\hhes}_0,\;\;\|\Xi_0\|=1.
\]
The integer $\nu$ is called {\em the order} of the pair  $(\lambda(t),
\Xi(t))$. A pair is called {\em degenerate} if its order is $>1$ and {\em nondegenerate} if it has order
$1$.Set
\[
Z^*=\{ (\lambda(t), \Xi(t))\in Z\; ;\; \lambda(t) \not\equiv 0\}.
\]
The complement $Z\setminus Z^*$ is  determined (according to Step 1) by $\ker
{\hhes}_t$ ($t>0$) and thus
\[
\#Z^*=\dim \ker {\hhes}_0 -\dim \ker{\hhes}_1.
\]
The spectral flow  $SF_+([\co])$  is then determined by
\be
SF_+([\co])=-\# \{(\lambda(t), \Xi(t))\in Z^*\; ;\; \lambda_\nu <0\}.
\label{eq: jet1}
\ee
To determine this integer we will distinguish two cases.

\medskip

\noindent {\bf The nondegenerate case} $\nu=1$. The equation ${\hhes}_t
\Xi(t)=\lambda(t) \Xi(t)$ implies
\[
{\hhes}_0\Xi_0=0\;\;\;{\hhes}_0\Xi_1+{\cal P}\Xi_0= \lambda_1 \Xi_0.
\]
This shows that $\lambda_1$ is a {\em nonzero} eigenvalue of the resonance
matrix ${\cal R}$ and moreover
\be
{\rm sign}\,\lambda_1 ={\rm sign}\, \lan {\cal R}\Xi_0, \Xi_0 \ran.
\label{eq: jet2}
\ee
In particular, the contribution to the spectral flow of the nondegenerate pairs
is equal to the number of negative eigenvalues of the resonance matrix ${\cal
R}$.   Thus we need to better understand the structure of the resonance form
\[
Q(\Xi)=\lan {\cal R}\Xi, \Xi\ran,\;\;\Xi\in \ker {\hhes}_0.
\]
Any $\Xi\in \ker {\hhes}_0$ decomposes as
\[
\Xi =\dps \oplus \ii \da \oplus \ii f
\]
where $\dps =\dps_-\oplus \dps_+\in \ker D_A$, $\da \in \Omega^1(N)$ is
harmonic  and $f$ is constant. All these objects are pulled back from the base
and  moreover 

\noindent $\bullet$ $\dot{\bar{\psi}}_-\in H^0(K^{1/2}-L)$, $\dps_+ \in
H^0(K^{1/2}+L)$.

\noindent $\bullet$ $\ii \da=\frac{1}{2}(\omega -\bar{\omega})$, $\omega\in H^0(K)$.

With these observations in place we have the following result.

\begin{lemma}
\[
Q(\Xi) =2^{1/2}f{\im}\lan \phi_+, \dps_+\ran - {\re}(\dbps_-\phi_+ \bar{\omega}).
\]
\label{lemma: jet1}
\end{lemma}
The proof of this lemma can be found in  Appendix D, equality (\ref{eq: d4}).

We see that (up to the positive factor $2^{1/2}$) the resonance form  is the direct sum of

\noindent (a)  a quadratic form $Q_1$ on ${\bR}\oplus H^0(K^{1/2}+L)$
\[
Q_1(f \oplus \dps_+) = f {\im}\lan \phi_+, \dps_+\ran.
\]
(b) a quadratic form $Q_2$ on $H^0(K^{1/2}-L)\oplus H^0(K)$ defined by
\[
Q_2(\dbps_- \oplus \omega) =-{\re}(\dbps_- \phi_+ \bar{\omega}).
\]
If we denote by $\dim_\pm$ the dimension of the positive/negative eigenspace of
a quadratic form  then
\[
\dim_- Q=\dim_-Q_1 +\dim_-Q_2.
\]

\medskip

\noindent{\sf The negative eigenspace of } $Q_1$. Set $V=H^0(K^{1/2}+L)$ and
$e_1=\phi_+$. Then
\[
\Omega(u,v) = {\im} \lan u, v \ran ,\;\;u,v\in V
\]
is a symplectic form on $V$. $Q_1$ is the quadratic form on ${\bR}\oplus V$
defined by
\[
Q_1(f\oplus v) =f \Omega(e_1, v)
\]
To determine its negative eigenspace  extend $e_1$ to a symplectic basis $e_1,
e_2, \ldots e_{2d-1}, e_{2d}$ where $d=\dim_{\bf C}V$. If $v=\sum_jv_je_j$ then
\[
Q_1(f\oplus v)=fv_2.
\]
This can be easily diagonalized and we get
\be
\dim_- Q_1 =1=\dim_+Q_1
\label{eq: jet3}
\ee
and
\be
\dim \ker Q_1= 2\dim_{\bf C}V -1 =2h_0(K^{1/2}+L)-1=2h_{1/2}(L)-1.
\label{eq: jet4}
\ee

\medskip

\noindent {\sf The negative eigenspace of} $Q_2$. Consider the multiplication
map
\[
{\bf m}: H^0(K^{1/2}-L)\ra H^0(K),\;\; \dbps_- \mapsto\dbps_- \phi_+.
\]
${\bf m}$ is obviously injective. (The implied inequality $\dim
H^0(K^{1/2}-L)\leq \dim H^0(K)$ is also a consequence of the classical Clifford
theorem.) Set $V=H^0(K)$ and $U={\rm Range}\, {\bf m}$.  $Q_2$ can be rewritten
as
\[
Q_2(\dbps, \omega)=-{\re}\lan {\bf m}\dbps, \omega\ran
\]
and thus it can be regarded as the quadratic form on $U\oplus V$
\[
Q_2(u\oplus v) =-{\re}\lan u, v\ran.
\]
This can be again easily diagonalized and  leads to the equalities
\be
\dim_- Q_2 = \dim_{\bf R}U= 2h_{1/2}(L^*)=\dim_+Q_2
\label{eq: jet5}
\ee
\be
\dim \ker Q_2= \dim_{\bf R}H^0(K) -\dim_{\bf R}U= 2g -2h_{1/2}(L^*).
\label{eq: jet6}
\ee
Summarizing, we deduce the following.

\noindent  {\bf A.}The spectral flow contribution of the nondegenerate pairs in $Z^*$ is
\[
-1-2h_{1/2}(-L).
\]
{\bf B.} The number of degenerate pairs $(\lambda(t), \Xi(t))\in Z$ is equal to
\[
\dim_{\bf R} \ker {\cal R}= 2h_{1/2}(L) +2g-2h_{1/2}(-L)-1=2(g-|n|)-1.
\]
Recall that $\dim_{\bf R} \ker {\hhes}_t =2(g-1-|n|)$ (if $t>0$)  and the pairs $(\lambda(t), \Xi(t))$ spanning $\ker {\hhes}_t$ do not contribute to the spectral flow. Hence there can be at most $\dim_{\bf R}\ker
{\cal R}-\dim \ker_{\bf R} {\hhes}_t=1$  degenerate pairs contributing to the spectral flow.

\noindent {\bf The degenerate case} $\nu > 1$. Set $d=\dim \ker {\hhes}_1$. We have  $d+1$  degenerate pairs  
\[
\{ (\lambda^k(t), \Xi^k(t))\; ;\; k=0, \ldots , d\}
\]
where the labeling is such that   $\ker {\hhes}_t ={\rm span}_{k\geq 1}\, (\Xi^k(t))$. Thus we need to determine the contribution of the pair $(\lambda^0(t), \Xi^0(t))$ to the spectral flow. First we claim that this pair has order two. Two achieve this we argue by contradiction.

Set
\[
S=\{\Xi_0\oplus \Xi_1\; ;\; \Xi_0\in \ker {\cal R},\;\;{\hhes}_0\Xi_1 +{\cal P}\Xi_0 =0\}.
\]
Using the perturbation series
\[
\lambda^0(t)=\lambda^0_\nu t^\nu+\cdots \;\; \nu \geq 3
\]
\[
\lambda^k(t)\equiv 0,\;\;\forall k=1,\ldots, d
\]
\[
\Xi^k(t)=\Xi^k_0+\Xi^k_1t+\Xi_2^kt^2+\cdots ,\;\;\;k=0,\ldots d
\]
we deduce that for all  $k=0, \ldots ,d$ 
\be
\left\{
\begin{array}{rcl}
{\hhes}_0\Xi^k_0&=&0 \\
{\hhes}_0\Xi^k_1+{\cal P}\Xi^k_0 &=& 0 \\
{\hhes}_0 \Xi^k_2+{\cal P}\Xi_1^k & = & 0
\end{array}
\right.
\label{eq: jet7}
\ee
Thus $\Xi_0\oplus \Xi^k_1\in S$ $\forall k$.  Taking the inner product with $\Xi_0^j$ in the last inequality we get
\be
\lan {\cal P}\Xi_1^k, \Xi_0^j\ran =0 ,\;\;\forall j,k=0,\ldots, d.
\label{eq: star}
\ee
Now observe the following elementary fact. Given 
\[
(\Xi_0, \Xi_1),\; (\Xi_0,\Xi'_1),\;(U_0, U_1)\in S
\]
i.e. $\Xi_0\in \ker {\cal P}$ and ${\hhes}_0\Xi_1={\hhes}_0 \Xi'_1=-{\cal P}\Xi_0$ then
\[
\lan {\cal P}\Xi_1, U_0\ran =\lan {\cal P}\Xi'_1, U_0\ran.
\]
Indeed,
\[
\lan {\cal P}\Xi_1 -{\cal P}\Xi'_1, U_0\ran =\lan \Xi_1-\Xi_1', {\cal P}U_0\ran
=0
\]
since $\Xi_1-\Xi_1'\in \ker {\hhes}_0$ and ${\cal P}U_0\perp \ker {\hhes}_0$.
In other word, the quantity
\[
{\cal B}(\Xi_0\oplus, \Xi_1, U_0\oplus U_1) = \lan {\cal P}\Xi_1, U_0\ran
\]
{\em depends bilinearly only upon $\Xi_0$ and $U_0$.} Thus it defines a
bilinear form on $\ker {\cal R}$ and one can check it is also symmetric. The
equality (\ref{eq: star}) implies
\[
{\cal B}(\Xi^j_0,\Xi^k_0)=0\;\;\forall j, k=0,\ldots, d
\]
i.e. ${\cal B}$  is trivial on $\ker {\cal R}$. We will show that this is not
the case thus establishing that the order of $\lambda^0(t)$ must be $2$.

Let  
\[
\Xi_0= \left[
\begin{array}{c}
0\oplus \phi_+ \\
0\\
0
\end{array}
\right]\in \ker {\hhes}_0.
\]
Using the identity (\ref{eq: d3}) in Appendix D we get
\[
{\cal P}\Xi_0=\left[
\begin{array}{c}
0\\
\ii|\phi_+|^2 \vfi\\
0
\end{array}\right] \in \left(\ker {\hhes}_0\right)^\perp.
\]
where $\vfi$ is the global angular form on $N$. Hence $\Xi_0\in \ker {\cal R}$. We have to solve
\[
{\hhes}_0\Xi_1+{\cal P}\Xi_0=0.
\]
If we write $\Xi_1=\dps \oplus \ii \da \oplus \ii f$ then the above equation
can be rewritten as
\be
\left\{
\begin{array}{rcl}
D_A\dps&=&0 \\
-\ast d \da +df +|\phi_+|^2\vfi&=&0 \\
d^*\da&=&0
\end{array}
\right.
\label{eq: jet9}
\ee
One can  say quite a lot about $\Xi_1$.  First note that since $|\phi_+|^2\vfi$ is co-closed (Appendix D) and $f \perp \{ constants\}$  we conclude that $f\equiv 0$.  We deduce
\be
\left\{
\begin{array}{rcl}
\ast d \da  &=&|\phi_+|^2\vfi\\
d^*\da & = &0
\end{array}
\right.
\ee
The above equation has an unique solution $\da$ orthogonal to the space of harmonic 1-forms. It is  given explicitly by
\be
\da=-\frac{1}{2\ell}|\phi_+|^2\vfi.
\label{eq: explicit}
\ee
Taking the inner product with $\da$ of the second equation of (\ref{eq: jet9})  we deduce that
\[
{\cal B}(\Xi_0) = \lan \Xi_1, {\cal P}\Xi_0\ran= -\lan {\hhes}_0\Xi_1, \Xi_1\ran
\]
\[
 =\int_N \lan \ast d\da, \da \ran dv_N =\int_N|\phi_+|^2\lan \vfi,
\da\ran dv_N=-\frac{1}{2\ell}\int_N |\phi_+|^2 dv_N.
\]
We conclude from (\ref{eq: explicit}) that
 \be
 {\cal B}(\Xi_0)= -{\rm sign}\, \ell.
 \label{eq: jet10}
 \ee
 This shows  that ${\cal B}$ is nontrivial. Since ${\cal B}$ can have at most
 one nonzero eigenvalue the above equality shows  that this eigenvalue has the
 same sign as $-\ell$.  If we now use the perturbation  equations for
 $\lambda^0(t)$ we obtain
 \[
 \left\{
 \begin{array}{rcl}
 {\hhes}_0\Xi^0_0&=&0 \\
 {\hhes}_0\Xi^0_1 +{\cal P}\Xi^0_0&=&0\\
 {\hhes}_0\Xi^0_2+{\cal P}\Xi^0_1 &=&\lambda^0_2\Xi^0_0
 \end{array}
 \right.
 \]
 We deduce 
 \[
 {\rm sign}\, \lambda^0_2={\rm sign}\,\lan {\cal P}\Xi^0_1, \Xi^0_0\ran ={\rm
 sign}\, {\cal B}(\Xi^0_0)=-{\rm sign}\, \ell.
 \]
 Thus the degenerate part contributes to the spectral flow only when $\ell >0$.

 We can now assemble all the informations we have collected so far in the
 following result.
 
 \begin{theorem}
 \[
 SF_+([\co])=\left\{
 \begin{array}{rlc}
 -2-2h_{1/2}(L^*)&,& {\rm if}\; \ell >0 \\
 -1-2h_{1/2}(L^*)&,& {\rm if}\; \ell <0
 \end{array}
 \right.
 \]
 \[
 SF_-([\co])=\left\{
 \begin{array}{rlc}
 -2-2h_{1/2}(L^*)&,& {\rm if}\; \ell <0 \\
 -1-2h_{1/2}(L^*)&,& {\rm if}\; \ell >0
 \end{array}
 \right.
 \]
 \label{th: jet}
 \end{theorem}
 
 \begin{remark}{\rm We have considered only the case $[\co]\in {\modu}_{\kappa,
 n}$ with $n<0$. The case $n>0$   can be approached by entirely similar
 methods and can be safely left to the reader. The corresponding formul{\ae} can be obtained from the above  by making the Serre duality substitution $L\ra L^*$.}
 \end{remark}

\subsection{Virtual dimensions}
In this final subsection we will show how one can use  Theorem \ref{th: jet} to compute
virtual dimensions of finite energy moduli spaces.  We will rely heavily on the
techniques of [MMR].

\begin{figure}
\centerline{\psfig{figure=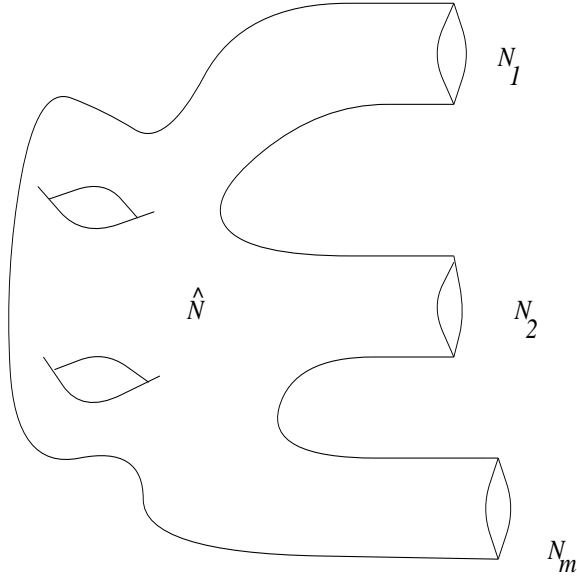,height=3in,width=3in}}
\caption{\sl Multiple cylindrical ends}
\label{fig: eta5}
\end{figure}

Consider a  $4$ manifold  $\hat{N}$ with a cylindrical end  isometric to
$[0,\infty)\times N$ where $N$ is  disjoint union of  nontrivial circle bundles
$\{N_j\; ;\; j=1,\ldots, m\}$ of degrees $\ell_j$ over Riemann surfaces
$\Sigma_j$ of genera $g_j\geq 1$ (see Figure \ref{fig: eta5}. 
Fix $spin$ structures on each of the Riemann
surfaces $\Sigma_j$ which induce by pullback $spin$ structures on $N_j$.  Next
fix a $spin^c$ structure $\hat{\si}$ on $\hat{N}$ which induces $spin^c$
structures $\sigma_j$ on $N_j$. Set 
\[
\sigma =\prod\sigma_j.
\]
The metric and compatible connections on the end of $\hat{N}$ are prescribed as
in dicated in \S 3.2.   This means  that as background connection on $N$ we use
the adiabatic connection $\nabla^0$. Consider  a finite energy solution 
$\hco=(\hat{\phi}, \hat{A})$ of the Seiberg-Witten equations  associated to the
structure $\hat{\si}$ and we assume  that along the neck it is in temporal gauge 
\[
\hco=\{t\mapsto\co(t)=(\phi(t),A(t))\}.
\]

The techniques of [MMR] work with no essential changes in the Seiberg-Witten
context and   show that  $[\co(t)]$ converges to $[\co_\infty] \in
\modu_\sigma$, where by $\modu_\sigma$ we denoted the Seiberg-Witten moduli space 
determined by the $spin^c$ structure $\si$ on $N$. The first conclusion we draw
from this fact is that $\sigma$ must be a pulled back $spin^c$ structure since
otherwise $\modu_\si =\emptyset$. Suppose this is indeed the case.

The moduli space $\modu_\si$ is a disjoint union
\[
\modu_\si =\prod \modu_{\si_j}
\]
and thus the asymptotic limit is a collection
\[
[\co_\infty]=([\co_1], \cdots, ,[\co_m]).
\]
Assume first that all the configurations $\co_j$ are irreducible 
\[
\co_j\in \modu_{\kappa_j,n_j}(\sigma_j).
\]
Again, to reduce the accounting job we consider that $n_j <0$ $\forall j$.  The
irreducibility condition implies that  $[\co(t)]$ converges {\em exponentially}
to its asymptotic limit.

We are interested in  describing a neighborhood of $\hco$ in the moduli space of
finite energy solutions and we will begin as in [MMR] by  studying a simpler
problem. 

Define $\hmodu([\co_\infty])$ the moduli space of  finite energy solutions 
with asymptotic limit $[\co_\infty]$. We want to  understand the structure of a
small neighborhood of $\hco \in \hmodu([\co_\infty])$.  More precisely, we would like
to compute the virtual dimension of such a neighborhood. This is achieved  using Kuranishi's deformation
picture of the moduli space which  requires  a suitable  functional framework.
Since the convergence to the asymptotic limit  is exponential one can use the
very convenient   weighted Sobolev spaces $L^{k,p}_w$ where $w$ is a very small
positive number.   The resulting  deformation complex can be described as in
Chap.8 of [MMR] and is
\be
0\ra
X_0\stackrel{ \hat{\Lie}_{\hco} }{\longrightarrow} X_1 \stackrel{\underline{sw} }{\longrightarrow}
X_2 \ra 0.
\label{eq: ku1}
\ee
where $X_0$ is the Lie algebra of the group of gauge  transformations on
$\hat{N}$  exponentially  converging to   1 along the neck
\[
X_0=L_w^{3,2}(\ii\Lambda^0T^*\hat{N})
\]
$X_1$ is the tangent space to the space of configurations of the
$4$-dimensional equations
\[
X_1=L_w^{2,2}(\hat{\bS}^+_{\hat{\si}}\oplus\ii\Lambda^1T^*\hat{N})
\]
$X_2$ is the space of ``obstructions''
\[
X_2=L_w^{1,2}(\hat{\bS}^-_{\hat{\si}}\oplus \ii \Lambda^2_+T^*\hat{N})
\]
$\hat{\Lie}=\hat{\Lie}_{\hco}$ is the infinitesimal  gauge group action at
$\hco$ and $\underline{sw}$ is the linearization  at $\hco$ of the $SW$-equations on
$\hat{N}$.

We can now form the  operator
\[
\hat{\cal O}_w:X_1\ra X_2\oplus X_0,\;\;\hat{\cal O}_w =\underline{sw}\oplus
\hat{\Lie}^{\ast_w}
\]
where ${}^{\ast_w}$ denotes the $L^2_w$-adjoint of $\hat{\Lie}$. This is an
elliptic operator and a computation {\em \`{a} la} [MMR] (Chap. 8) shows that along the neck 
it has the $APS$ form
\[
\hat{\cal O}_w=something \times(\nabla_t-{\cal O}_w) 
\]
where
\[
{\cal O}_w\left[
\begin{array}{c}
{\dps}\oplus{\ida} \\
{\bf i}f 
\end{array}
\right ] = \left[
\begin{array}{clcl}
D_A\dps  &                               &+ &  {\bc}(\ida)\phi-\ii f\phi   \\
                 &-\ii\ast  d\da +\ii df & +&\dta(\phi, \dps) \\
                  &\ii d^*\da -2w\ii f& + & \ii \im \lan \phi, \dps\ran  
\end{array}
\right]
\]
and $[\phi, A]=[\co_\infty]$.  Note that ${\hhes}_1(\co_\infty)={\cal
O}_w\!\mid_{w=0}$. $\hat{\cal O}_w$ is a Fredholm operator and its index (over
${\bR}$) is 
equal to the virtual dimension  of a small neighborhood of $[\hco] $
in $\hmodu([\co_\infty])$.  Remark \ref{rem: 110} at the end of \S 1.3 shows that  the index of $\hat{\cal
O}_w$ is equal to the $APS$ index of $\hat{\cal O}_w$.

Denote by ${\cal A}$ the   anti-selfduality operator on $\hat{N}$
\[
{\cal A}=d_+\oplus d^*:\ii\Omega^1(\hat{N})\ra \ii\Omega^2(\hat{N})\oplus
\ii\Omega^0(\hat{N}).
\]
Using the connection $\hat{A}$  and the $spin^c$ structure $\hsi$ we can form
the Dirac operator 
\[
\hat{D}_{\hat{A}}:\Gamma(\hat{\bS}^+_{\hat{\si}})\ra
\Gamma(\hat{\bS}^-_{\hat{\si}}).
\]
Along the neck  the direct sum ${\cal N}_{\hco}=\hat{D}_{\hat{A}}\oplus {\cal A}$ has the
$APS$ form
\[
{\cal N}=something \times \left(\nabla_t -{\hhes}_0(\co_\infty)\right). 
\]
Using the excision formula (\ref{eq: ex2})  in Remark \ref{rem: 110} we deduce
\be
{\rm ind}_{APS}(\hat{\cal O}_w) = {\rm ind}_{APS}({\cal N})-
SF({\hhes}_0\ra {\hhes}_1) - SF({\hhes}_1\ra {\cal O}_w)
\label{eq: vd0}
\ee 
where  $SF(A\ra B)$ denotes $SF(A+t(B-A),\; t\in[0,1])$.  All the indices
and the spectral flows above are   {\em real} quantities.

We now proceed to determine the three terms in the right-hand side of the above
formula.   

Corollary \ref{cor: d1} shows that the third term above vanishes. The second spectral term can be rewritten as
\be
SF({\hhes}_0\ra {\hhes}_1)=\sum_j SF_+([\co_j])
\label{eq: vd1}
\ee
We denote by $\rho_{asd}$ (resp. $\rho_{dir}$) the index densities of ${\cal
A}$ (resp. $\hat{D}_{\hat{A}}$)
\[
\rho_{asd}= -\frac{1}{2}\left({\bf e}(\hat{N})+{\bf L}(\hat{N})\right).
\]
where ${\bf e}(\hat{N})$ and ${\bf L}(\hat{N})$ denote respectively the Euler
and the $L$-genus forms on $\hat{N}$ constructed using the Levi-Civita
connection. Also
\[
\rho_{dir} = 2 \hat{\bf A}(\hat{\nabla}^0) \wedge \exp(\frac{1}{2}c_1(\det \hat{\si}))
\]
where on $\det \hat{\si}$ we used the connection induced by $\hat{A}$.  The
factor 2 appears since we are interested in the {\em real} index of $\hat{D}$. 
The $\hat{\bf A}$-genus form is  computed using the metric compatible  connection
$\hat{\nabla}^0$
which along the neck has the product form $dt\otimes \partial_t+\nabla^0$. For simplicity, set
$c(\hat{A})=c_1(\det \hat{\si})$.

On  a $4$-manifold the above equality has a simpler form
\[
\rho_{dir}=\frac{1}{4}(c(\hat{A})^2 -{\bf L}(\hat{\nabla}^0)).
\]
The $\xi$ invariant of ${\hhes}_0$   is  the sum  $\xi({\cal
A}\!\mid_N)+2\xi(D_A)$ (the factor 2 is present  for reality reasons).
\[
\xi({\hhes}_0)=\frac{1}{2}( \dim_{\bf R}\ker {\hhes}_0 -\eta_{sign} +
2\eta(D_A))
\]
where $\eta_{sign}$ denotes the eta invariant of the odd signature operator. We deduce
\[
{\rm ind}_{APS}({\cal N})=\int_N (\rho_{asd}+\rho_{dir})-\xi({\hhes}_0)
\]
\[
=-\frac{1}{2}\int_{\hat{N}} {\bf e} -\frac{1}{2}\left(\int_{\hat{N}} {\bf L}
-\eta_{sign}\right)+\frac{1}{4}\int_{\hat{N}}(c^2(\hat{A})-{\bf L}({\nah}^0))
\]
\[
-\frac{1}{2}\dim_{\bf R}\ker{\hhes}_0 -\eta(D_A)
\]
\[
=-(\chi(\hat{N}) +{\rm sign}\,(\hat{N}))/2 +\frac{1}{4}\int_{\hat{N}}(c^2(\hat{A})-{\bf
L}({\nah}^0))
\]
\[
-\sum_j\dim_{\bf C}\ker D_{A_j}-\frac{1}{2}\sum_j(2g_j+1) -\frac{1}{6}\sum_j\ell_j
\]
Using (\ref{eq: vd0}) and (\ref{eq: vd1}) we deduce
\[
{\rm ind}\,(\hat{\cal O}_w)= -\frac{1}{2}(\chi(\hat{N}) +{\rm sign}\,(\hat{N}))
+\frac{1}{4}\int_{\hat{N}}(c^2(\hat{A})-{\bf
L}({\nah}^0))
\]
\be
-\sum_j\left(\dim_{\bf C}\ker D_{A_j}+SF_+([\co_j]) \right)
-\frac{1}{2}\sum_j(1+ 2g_j)-\sum_j\frac{\ell_j}{6}.
\label{eq: vd2}
\ee
This formula can be further simplified.   We can replace  the integral
of ${\bf L}(\nah^0)$ with the integral of ${\bf L}(\hat{N})$ plus a correction
term given by the  second transgression formula. Assume for simplicity that all
components of $N$ have fibers of the same radius $r$ and the bases have common
area $\pi$. We get
\be
\int_{\hat{N}}{\bf
L}(\nah^0)-\int_{\hat{N}}{\bf
L}(\hat{N})=\sum_j\frac{2\ell_j}{3}(\ell_j^2r^4-\chi_jr^2)
\label{eq: vd3}
\ee
where $\chi_j=2-2g_j$.   Denote  by $\eta_j$ the signature eta invariant of
$N_j$. This was computed in [Ko] and [O] and we  have
\[
\eta_j=-{\rm sign}\,(\ell_j) -\frac{2\ell_j}{3}(\ell_j^2r^4-\chi_jr^2)
+\frac{\ell_j}{3}.
\]
If we set $\eta=\sum\eta_j$ we deduce from (\ref{eq: vd3})
\[
\int_{\hat{N}}{\bf L}(\nah^0)+\eta-\int_{\hat{N}}{\bf L}(\hat{N})=\sum_j(\ell_j/3-{\rm sign}\,(\ell_j)).
\]
The term 
\[
\eta-\int_{\hat{N}}{\bf L}(\hat{N})
\]
is equal to $-{\rm sign}\,(\hat{N})$ so that we  deduce
\[
\int_{\hat{N}}{\bf L}(\nah^0)={\rm sign}\,(\hat{N}) +\sum_j(\frac{\ell_j}{3} -{\rm
sign}\,(\ell_j))).
\]
If we use this equality in (\ref{eq: vd2}) we deduce
\[
{\rm ind}\,(\hat{\cal O}_w)= 
\frac{1}{4}\int_{\hat{N}}c^2(\hat{A})-\frac{1}{4}(2\chi(\hat{N})+3 {\rm
sign}\,(\hat{N}) )
\]
\be
-\sum_j\left(\dim_{\bf C}\ker D_{A_j}+SF_+([\co_j]) \right)
-\frac{1}{2}\sum(2g_j+1)-\frac{1}{4}\sum_j(\ell_j -{\rm sign}\,(\ell_j)).
\label{eq: vd4}
\ee

To find the virtual dimension $\dim_v(\hco)$ of a neighborhood of $\hco$ in the {\em entire}
moduli space $\hmodu$  we only have to add the dimensions of the  asymptotic
limit sets $\dim {\modu}_{\kappa_j, n_j}=2(g_j -1+n_j)$ (recall that we have
assumed $n_j<0$).
\[
\dim_v(\hco)=
\frac{1}{4}\int_{\hat{N}}c^2(\hat{A})-\frac{1}{4}(2\chi(\hat{N})+3 {\rm
sign}\,(\hat{N}) )
\]
\[
-\sum_j\left( \dim_{\bf C}\ker D_{A_j}+SF_+([\co_j]) \right)
\]
\be
+\frac{1}{2}(2g_j-1)+2\sum_j(n_j-1)-\frac{1}{4}\sum_j(\ell_j -{\rm sign}\,(\ell_j)).
\label{eq: vd5}
\ee
We can now use Theorem \ref{th: jet} in the form
\[
SF_+(\co_j)=-1-2h_{1/2}(L_j^*)-{\ve}_j
\]
where 
\[
{\ve}_j=\frac{1}{2}( 1+{\rm sign}\,(\ell_j)).
\]
Since $\dim_{\bf C}\ker D_{A_j}=h_{1/2}(L_j)+h_{1/2}(L^*_j)$ we deduce
\[
\dim_{\bf C}\ker D_{A_j}+SF_+([\co_j])=h_{1/2}(L_j)-h_{1/2}(L^*_j)-1-{\ve}_j=n_j-1-{\ve}_j
\]
Using this in (\ref{eq: vd5}) we get
\[
\dim_v(\hco)=\frac{1}{4}\int_{\hat{N}}c^2(\hat{A})-\frac{1}{4}(2\chi(\hat{N})+3 {\rm
sign}\,(\hat{N}) )
\]
\be
+\sum_j({\ve}_j+n_j-1)+\frac{1}{2}\sum_j(2g_j-1)-\frac{1}{4}\sum_j (\ell_j -{\rm sign}\,(\ell_j))
\label{eq: vd6}
\ee
If we define  the boundary contribution of  the asymptotic limit $[\co_j]$ by
\[
\beta([\co_j])=({\ve}_j+n_j-1) +\frac{1}{2}(2g_j-1)-\frac{1}{4}(\ell_j -{\rm sign}\,(\ell_j))
\]
then we can rewrite  formula (\ref{eq: vd6}) as
\be
\dim_v(\hco)=\frac{1}{4}\int_{\hat{N}}c^2(\hat{A})-\frac{1}{4}(2\chi(\hat{N})+3 {\rm
sign}\,(\hat{N}) ) +\sum_j\beta([\co_j]).
\label{eq: vd11}
\ee
The first two terms in (\ref{eq: vd11}) represent  formally the expression which computes
the virtual dimension of the Seiberg-Witten moduli spaces on a closed compact
$4$-manifold.

We can now apply this formula to the special situation of tunnelings. In this
case $\hat{N}={\bR}\times N$ and thus $m=2$, $g_1=g_2=g$, $\ell_1+\ell_2=0$.
Moreover ${\rm sign}\,(\hat{N})=\chi(\hat{N})=0$.  The equality (\ref{eq: vd6})    
gives the virtual dimension  of the space  of tunnelings between
$\modu_{\kappa, n_1}$ (at $-\infty$) and ${\modu}_{\kappa, n_2}$ (at $+\infty$).  Denote this dimension by
$\tau(\kappa; n_1,n_2)$. We have
\[
\tau(\kappa; n_1,n_2)=\int_{\hat{N}}c(\hat{\si})^2  +n_1+n_2 +2g-2.
\]
The integral term can be computed via transgression exactly as in the third
transgression formula.  It suffices to pick arbitrary $[\co_j]=[\phi_j,A_j]\in \modu_{\kappa, n_j}$, $j=1,2$. and an arbitrary path of connections $A(t)$ such  that $A(0)=A_1$ and $A(1)=A_2$. We choose $A_j$ such that
\[
\frac{\ii}{2\pi}F_{A_j}=\frac{n_j}{\pi}dv_\Sigma
\]
and
\[
A_2-A_1=\ii c\vfi
\]
where $c=(n_2-n_1)/\ell$.  $A(t)$ will be the affine path $A_1+t\ii c\vfi$.  
(The connection on $\det \hat{\si}$ will be $\hat{A}^{\otimes 2}$ since $\det
{\bS}_L=L^2)$).  We get a connection $\hat{A}$ on $[0,1]\times N$. 
 A computation  entirely similar to the one in the third transgression formula 
leads to the equality
\[
\frac{1}{2}\int_{\hat{N}}c(\hat{\si})^2 =-\frac{1}{4\pi^2}\int_{[0,1]\times
N}F_{\hat{A}}=\frac{n_1^2-n_2^2}{\ell}.
\]
Thus we get
\[
\tau(\kappa ;n_1,n_2)=\frac{n_1^2-n_2^2}{\ell}+n_1+n_2+2g-2.
\]
This agrees with Corollary 1.0.5 of [MOY]. (In the notation of [MOY],
$e_j=n_j+g-1$.)

We have omitted from our discussion the case when one  (or several) asymptotic
limits $[\co_j]$ is reducible.  The degenerate case, $g_j-1\equiv 0$ (mod
$\ell_j$), requires special care and will not be discussed here. In the remaining cases  the problem is  actually
simpler than the case dealing with irreducible limits.  

First of all  the convergence to such a  nondegenerate reducible continues to be exponential and
thus we can use the same   functional framework as above. Assume for simplicity
the boundary has only one component. We have to compute
the $APS$ index  of a new operator $\hat{\cal O}_w$ which  along the neck has
the form
\[
\hat{\cal O}_w=something \times (\nabla_t-{\cal O}_w)
\]
where this time
\[
{\cal O}_w\left[
\begin{array}{c}
{\dps}\oplus{\ida} \\
{\bf i}f 
\end{array}
\right ] = \left[
\begin{array}{cl}
D_A\dps  &                                  \\
                 &-\ii\ast  d\da +\ii df  \\
                  &\ii d^*\da -2w\ii f   
\end{array}
\right]
\]
(The spinor part $\phi$ of the asymptotic limit $[\co]=[\phi, A]$ is zero and thus ${\cal
P}_\phi\equiv 0$.) Thus  
\[
{\rm ind}_{aps}(\hat{\cal O}_w)= {\rm ind}_{aps}({\cal N})- SF({\cal O}_0\ra{\cal O}_w).
\]
The spectral flow contribution is easy to determine.     The only eigenvalue of
${\cal O}_{tw}$  contributing to the spectral flow is $-2wt\!\mid_{t=0}$ with
a single eigenfunction $\dps \oplus \ii \da \oplus \ii f$, where $\dps=0$,
$\ii\da=0$ and $f\equiv 1$.  Hence
\[
{\rm ind}_{aps}(\hat{\cal O}_w)= {\rm ind}_{aps}({\cal N}) +1.
\]
The index of ${\cal N}$ can be determined as above using instead the eta invariant of the adiabatic operator coupled with the flat
connection $A$ determined in \S 2.4. The nondegeneracy condition also implies  $\ker D_A=0$. Hence this eta invariant  is twice the $\xi$-invariant in (\ref{eq:
eta3}) so that
\[
\eta(D_A)=\frac{\ell}{6}+\frac{\kappa^2}{\ell}-\kappa {\rm sign}\,(\ell).
\]
We deduce (using $\dim_{\bf R}\ker {\hhes}_0=b_0(N)+b_1(N)=2g+1$)
\[
{\rm ind}_{APS}({\cal O})_w= 1-(\chi(\hat{N}) +{\rm sign}\,(\hat{N}))/2
+\frac{1}{4}\int_{\hat{N}}(c^2(\hat{A})-{\bf
L}({\nah}^0))
\]
\[
 -\frac{1}{2}(1+ 2g)-\frac{\ell}{6} -\frac{\kappa^2}{\ell}+\kappa {\rm sign}\,(\ell).
\label{eq: vd7}
\]
Arguing as in the irreducible case we deduce
\[
{\rm ind}\,(\hat{\cal O}_w)= 
\frac{1}{4}\int_{\hat{N}}c^2(\hat{A})-\frac{1}{4}(2\chi(\hat{N})+3 {\rm
sign}\,(\hat{N}) )
\]
\[
 -\frac{1}{2}(2g-1)-\frac{1}{4}\left(\ell -{\rm
sign}\,(\ell)\right)-\frac{\kappa^2}{\ell}+
\kappa {\rm sign}\,(\ell).
\label{eq: vd12}
\]
To find the virtual dimension of a neighborhood in the entire  moduli space
${\hmodu}$ we proceed as in Sec. 8.5 of [MMR].  We  need to add the dimension of the reducible limit set (which is a $2g$-torus) and subtract the dimension of the stabilizer of the  asymptotic limit (which is $S^1$). We deduce
\[
\dim_v (\hco)= 
\frac{1}{4}\int_{\hat{N}}c^2(\hat{A})-\frac{1}{4}(2\chi(\hat{N})+3 {\rm
sign}\,(\hat{N}) )
\]
\be
 +\frac{1}{2}(2g-1)-\frac{1}{4}(\ell -{\rm sign}\,(\ell))-\frac{\kappa^2}{\ell}+\kappa {\rm sign}\,(\ell).
\label{eq: vd8}
\ee
We see that the boundary contribution of a reducible limit is
\be
\beta([\co])=\frac{2g-1}{2}-\frac{\ell-{\rm sign}\,(\ell)}{4} - \frac{\kappa^2}{\ell}+\kappa {\rm sign}\,(\ell).
\label{eq: vd9}
\ee
We can now easily write the virtual dimension $\tau(\kappa:0,n)$of the space of tunnelings from a
reducible solution $[\co_1]$ to an irreducible one $[\co_2]\in \modu{\kappa,
n}$.  It is 
\[
\tau(\kappa;0,n) =\int_N{\bf T}c_1^2(A_2,A_1) +\beta([\co_1])+\beta([\co_2]).
\]
where ${\bf T}$ stands for the transgression form. We denote by $\ell$ the degree of the boundary at $+\infty$. Then as in the
third transgression formula we get
\[
\int_N{\bf T}c_1^2(A_2,A_1)= -\frac{n^2}{\ell}.
\]
Also
\[
\beta([\co_2])=(\frac{1+{\rm sign}\,(\ell)}{2} +n-1)+\frac{2g-1}{2}-\frac{\ell -{\rm sign}\,(\ell)}{4}
\]
and (using the opposite orientation at $-\infty$)
\[
\beta([\co_1])=\frac{2g-1}{2}+\frac{\ell -{\rm
sign}\,(\ell)}{4}+\frac{\kappa^2}{\ell}-\kappa {\rm sign}\,(\ell).
\]
We get
\[
\tau(\kappa;0,n)=\frac{\kappa^2-n^2}{\ell} +2g-2 +n + \frac{1+{\rm sign}\,(\ell)}{2}
\]

\appendix 
\section{Proof of the first transgression formula}
\setcounter{equation}{0}
 Let $0 < \rho \ll r \ll 1$.  {\em The parameter $r$   will  stay fixed throughout this section.
 $\rho$  will eventually go to } $0$.   We have a {\em fixed}  local frame
 $(\zeta, \zeta_1, \zeta_2)$. This is {\em not an orthonormal frame} for
 either of the metric $h_r$ or $h_\rho$ but it is an orthogonal frame.   
 
 We first want to compute the 1-forms  associated  by the above frame to  the 
 connections $\nabla^r$ and $\nabla^\rho$.  We denote them  by $\Gamma_r$ and
 resp. $\Gamma_\rho$.  Using the equalities (\ref{eq: structural}) and
 (\ref{eq: str1}) we get after some simple manipulations that
 \[
 \Gamma_\rho=\left[
 \begin{array}{ccc}
 0 & \lambda \vfi^2 & -\lambda \vfi^1 \\
 -\lambda\vfi^2 & 0 & -r^2\lambda \vfi -\kappa \vfi^1 \\
 r^2\lambda \vfi^1 & r^2\lambda \vfi \kappa \vfi^1 & 0
 \end{array}
 \right]
 \]
 We get a similar result for $\Gamma_r$. Set $\Xi_\rho =\Gamma_\rho -\Gamma_r$
 and $\Xi_0=\lim_{\rho\ra 0}\Xi_\rho$. Note that
 \[
 \Xi_0=\left[
 \begin{array}{ccc}
 0 & 0 & 0 \\
 r^2\lambda \vfi^2 & 0 & r^2\lambda\vfi \\
 
 -r^2\lambda\vfi^1 & -r^2 \lambda\vfi & 0 
 \end{array}
 \right]
 \]
 We have  to compute $\lim_{\rho \ra 0}T\hat{\bf A} (\nabla^\rho, \nabla^r)$. Note that
 \[
 T\hat{\bf A}(\nabla^\rho, \nabla^r)=\frac{1}{96\pi^2}{\rm tr}\,\left\{
 \Xi_\rho \wedge \left(\Omega_r + \frac{1}{2}(d\Xi_\rho+ \Gamma_r\wedge\Xi_\rho
 +\Xi_\rho \wedge \Gamma_r) +\frac{1}{3}\Xi_\rho \wedge \Xi_\rho \wedge
 \Xi_\rho\right) \right\}.
 \]
 Above, $\Omega_r$ is the curvature $2$-form of the connection $\nabla^r$. By letting $\rho \ra 0$ in the above equality we deduce
 \[
 \lim_{\rho \ra 0}T\hat{\bf A} (\nabla^\rho, \nabla^r)
 \]
 \be
 =\frac{1}{96\pi^2}{\rm tr}\,\left\{
 \Xi_0 \wedge \left(\Omega_r + \frac{1}{2}(d\Xi_0+ \Gamma_r\wedge\Xi_0
 +\Xi_0 \wedge \Gamma_r) +\frac{1}{3}\Xi_0 \wedge \Xi_0 \wedge
 \Xi_0\right) \right\}.
 \label{eq: pont1}
 \ee
 We now proceed to describe  each of the above constituents, one by one. Since
 $\Omega_r =d\Gamma_r +\Gamma_r \wedge \Gamma_r$ we deduce  using (\ref{eq:
 di0}), (\ref{eq: di1}) and (\ref{eq: di2}) that
 \[
 \Omega_r=\left[
 \begin{array}{ccc}
 0& r^2\lambda^2 \vfi\wedge \vfi^1 & -\lambda^2 r^2 \vfi^2 \wedge \vfi \\
 -\lambda^2 r^4 \vfi \wedge \vfi^1 & 0 & -(\kappa^2 +3\lambda^2 r^2)\vfi^1\wedge
 \vfi^2\\
 \lambda^2r^4 \vfi^2\wedge\vfi & (\kappa^2 + 3\lambda^2 r^2) \vfi^1\wedge \vfi^2
 & 0 
 \end{array}
 \right]
 \]
 Then
 \be
 \Xi_0 \wedge \Omega_r =\left[
 \begin{array}{ccc}
 0 &\ast & \ast \\
 \ast & (4\lambda^3r^4 +\lambda\kappa^2 r^2) & \ast \\
 \ast & \ast & (4\lambda^3r^4 +\lambda\kappa^2 r^2)
 \end{array}
 \right]\vfi\wedge \vfi^1 \wedge \vfi^2.
 \label{eq: tr1}
 \ee
 \[
 d\Xi_0=\left[
 \begin{array}{ccc}
 0&0&0 \\
 0&0& 2\lambda^2r^2 \\
 -\lambda\kappa r^2 & -2\lambda^2r^2 & 0
 \end{array}\right]\vfi^1\wedge \vfi^2.
 \]
 Then
 \be
 \Xi_0\wedge d\Xi_0=\left[
 \begin{array}{ccc}
 0& \ast &\ast\\
 \ast & -2\lambda^3r^4 & \ast \\
 \ast & \ast & -2\lambda^3 r^4 
 \end{array}
 \right]\vfi\wedge \vfi^1\wedge \vfi^2.
 \label{eq: tr2}
 \ee
 Simple manipulations yield
 \be
 \Xi_0\wedge \Xi_0 \wedge \Gamma_r +\Xi_0\wedge \Gamma_r\wedge \Xi_0=\left[
 \begin{array}{ccc}
 0 & \ast & \ast \\
 \ast & -2\lambda^3r^4 & \ast \\
 \ast & \ast & -2\lambda^3r^4
 \end{array}
 \right]\vfi\wedge\vfi^1\wedge \vfi^2.
 \label{eq: tr3}
 \ee
 An immediate computation (eased by the large number of vanishing entries in
 $\Xi_0$) shows that ${\rm tr}\, (\Xi_0\wedge \Xi_0 \wedge \Xi_0)=0$.  By combining (\ref{eq: pont1}) with
 (\ref{eq: tr1})-(\ref{eq: tr3}) we get
 \[
 \lim_{\rho \ra 0}T\hat{\bf A} (\nabla^\rho,
 \nabla^r)=\frac{1}{96\pi^2}(4\lambda^3r^4+2\lambda\kappa^2r^2)\vfi \wedge
 \vfi^1\wedge \vfi^2.
 \]
 The first transgression formula now follows by integrating over $N$ and using
 the equalities $\lambda=-\ell$, $\kappa^2=4(g-1)$ and
 \[
 \int_N\vfi\wedge \vfi^1\wedge \vfi^2 =2\pi^2.
 \]

\section{Proof of the second transgression formula}
\setcounter{equation}{0}
We have to compute $\lim_{t\ra 0}T\hat{\bf A}(\nabla^{r,t}, \nabla^r)$. As in
Appendix  A we get
\[
\lim_{t\ra 0}T\hat{\bf A}(\nabla^{r,t}, \nabla^r)
\]
\be
=\frac{1}{96\pi^2}{\rm tr}\,\left\{
 \Xi_0 \wedge \left(\Omega_r + \frac{1}{2}(d\Xi_0+ \omega_r\wedge\Xi_0
 +\Xi_0 \wedge \omega_r) +\frac{1}{3}\Xi_0 \wedge \Xi_0 \wedge
 \Xi_0\right) \right\}
 \label{eq: pont2}
 \ee
 where $\Omega_r=d\omega_r +\omega_r\wedge \omega_r$ and (using (\ref{eq:
 str1}) and (\ref{eq: str2}))
 \[
 \Xi_0=\omega_{r,0}-\omega_r =-\lambda_r\left[
 \begin{array}{ccc}
 0 & \vfi^2 &-\vfi^1 \\
 -\vfi^2 & 0 & -\vfi_r \\
 \vfi^1 & \vfi_r & 0
 \end{array}
 \right].
 \]
 Using (\ref{eq: di0}), (\ref{eq: di1}) and (\ref{eq: di2}) we get after some
 simple manipulations
 \[
 \Omega_r=\left[
 \begin{array}{ccc}
 0 & \lambda_r^2 \vfi_r\wedge \vfi^1 & -\lambda_r^2 \vfi^2\wedge \vfi_r \\
 -\lambda_r^2 \vfi_r\wedge \vfi & 0 & -(3\lambda_r^2 +\kappa^2)\vfi^1\wedge
 \vfi^2 \\
 \lambda_r^2 \vfi^2\wedge \vfi_r &(3\lambda_r^2 +\kappa^2)\vfi^1\wedge
 \vfi^2 & 0 
 \end{array}
 \right].
 \]
 Then
 \be
 \Xi_0\wedge \Omega_r =-\lambda_r\left[
 \begin{array}{ccc}
 -2\lambda_r^2 & \ast & \ast \\
 \ast & -(4\lambda_r^2 +\kappa^2) & \ast \\
 \ast & \ast & -(4\lambda_r^2 +\kappa^2)
 \end{array}
 \right]\vfi_r\wedge \vfi^1\wedge \vfi^2.
 \label{eq: tr5}
 \ee
 \[
 d\Xi_0=-\lambda_r\left[
 \begin{array}{ccc}
 0 & 0 & -\kappa \\
 0 & 0 & -2\lambda_r \\
 \kappa & 2\lambda_r & 0 
 \end{array}
 \right]\vfi^1\wedge \vfi^2.
 \]
 \be
 \Xi_0 \wedge d\Xi_0 =\lambda_r^2 \left[
 \begin{array}{ccc}
 0 & \ast & \ast \\
 \ast & -2\lambda_r & \ast \\
 \ast & \ast & -2\lambda_r
 \end{array}
 \right]\vfi_r\wedge \vfi^1\wedge \vfi^2.
 \label{eq: tr6}
 \ee
 \be
 \Xi_0 \wedge \Xi_0 \wedge \Xi_0 =2\lambda_r^3\left[
 \begin{array}{ccc}
 1 & \ast & \ast \\
 \ast & 1 & \ast \\
 \ast & \ast & 1
 \end{array}
 \right]\vfi_r \wedge \vfi^1 \wedge \vfi^2.
 \label{eq: tr7}
 \ee
 \be
 \Xi_0\wedge\omega_r \wedge \Xi_0 +\Xi_0 \wedge \Xi_0 \wedge
 \omega_r=-4\lambda_r^3\left[
 \begin{array}{ccc}
 1 & \ast & \ast \\
 \ast & 1 & \ast \\
 \ast & \ast & 1
 \end{array}
 \right]\vfi_r \wedge \vfi^1 \wedge \vfi^2.
 \label{eq: tr8}
 \ee
 Putting together (\ref{eq: tr5})-(\ref{eq: tr8}) we deduce
 \[
 \lim_{t\ra 0}T\hat{\bf A}(\nabla^{r,t},
 \nabla^r)=\frac{1}{96\pi^2}(4\lambda_r^3+2\kappa^2\lambda_r)\vfi_r\wedge \vfi^1\wedge
 \vfi^2
 \]
 \[
 =\frac{1}{96\pi^2}(4\lambda^3r^4 +2\lambda \kappa^2 r^2)\vfi \wedge \vfi^1 \wedge
 \vfi^2.
 \]
 We now conclude exactly as in Appendix A.  $\Box$

\section{Elementary computation of the eta invariants}
\setcounter{equation}{0}
We include here an elementary derivation  of the equality (\ref{eq: adiaeta}).  For brevity we  present a proof  only for circle bundles over  smooth Riemann surfaces but the arguments extend  to the more general case of Seifert manifolds, i.e. smooth circle bundles over Riemann V-surfaces  (2-orbifolds). The changes  from the smooth to the orbifold case are only cosmetic (``orbify'' everything i.e. add the prefix $V$ to all the intervening geometric objects and use  known $V$-theorems: V-Riemann-Roch, V-Serre duality etc.).  

Our circle bundle $N$ equipped with the metric described in $\S 2.1$  determines a hermitian metric and compatible connection on  a degree $\ell$  hermitian line bundle over $\Sigma$. This connection determines a holomorphic structure and we denote by $L_0$  the holomorphic  line bundle thus obtained.  Consider another line bundle $L\ra \Sigma$  of degree $k$ equipped with a hermitian metric and compatible connection $B$.   The connection $B$ determines a holomorphic structure on $L$  and    we  will denote  by  $h(L)$ the  dimension of the space of holomorphic sections.     

Consider now the $spin^c$ structure on $N$ whose associated spinor bundle is 
\be
{\bS}_k ={\can}^{-1}\otimes \pi^*L \oplus \pi^* L.
\label{eq: c0}
\ee
In terms of the notations in $\S 2.4$ we have ${\bS}_k ={\bS}\otimes {\can}^{-1/2}\otimes L$.  Note that ${\bS}_k$ makes no reference to a choice of  $spin$ structure on the base $\Sigma$ and thus  a similar object can also be defined  when $N$ is a Seifert fibration over a not necessarily $spin$-orbifold.  This is very similar to  the case of $spin^c$ structures over (even dimensional) almost complex manifolds. 

Using the pullback of $B$  on $\pi^*L$, the pullback of the Levi-Civita connection on $K_\Sigma$ to ${\can}$ and the adiabatic Levi-Civita connection $\lim_{t\ra 0}\nabla^{r,t}$ we obtain as in $\S 2.4$ the adiabatic Dirac operator  $D=D_B$.   For each $\mu\in {\bR}$ define
\[
V_\mu =\ker (\mu-D), \;\;v_\mu=\dim V_\mu.
\]
We want to compute the eta function  of $D$
\be
\eta(s)=\eta_D(s)=\sum_{\mu>0}\frac{v_\mu-v_{-\mu}}{\mu^s}.
\label{eq: c1}
\ee
Recall  ($\S 2.4$) that $D$  has a decomposition $D=Z+T$. Now define
\[
E_\mu=\{v\in V_\mu\; ;\; ZTv=0\},\;\; e_\mu=\dim E_\mu.
\]
Denote by $E_\mu^\perp$ the orthogonal complement of $E_\mu$ in $V_\mu$.  Since $\{Z,T\}=0$  we also have $\{D, ZT\}=0$ and it is  easy to check that $ZT(E_\mu^\perp )\subset E_{-\mu}^\perp$. The definition of $E_\mu$ implies that the induced map $ZT:E_\mu^\perp \ra E_{-\mu}^\perp$ is injective. Thus $\dim E^\perp_\mu\leq \dim E_{-\mu}^\perp$ and by symmetry $\dim E^\perp_\mu = \dim E_{-\mu}^\perp$. Using this in (\ref{eq: c1}) we deduce
\be
\eta(s) =\sum_{\mu>0}\frac{e_\mu-e_{-\mu}}{\mu^s}.
\label{eq: c2}
\ee
After some elementary manipulation which can be safely left to the reader we   
deduce that $E_\mu$ is both $Z$ and $T$ invariant. Since $\{Z,T\}=ZT=0$ we deduce
 that  $Z$, $T$ {\em commute} as operators on $E_\mu$ and moreover, $Z+T=D\equiv \mu$ on $E_\mu$.  Standard spectral 
 theory for commuting  symmetric operators  implies  that $E_\mu$  admits an orthogonal decomposition $E_\mu=F_\mu \oplus B_\mu$ with respect to which $Z$ and $T$ have the block decompositions
\[
Z=\left[
\begin{array}{cc}
\mu & 0 \\
0 & 0 
\end{array}
\right],\;\;T=\left[
\begin{array}{cc}0 & 0\\
0 & \mu
\end{array}
\right].
\]
Set  $f_\mu=\dim F_\mu$ and $b_\mu=\dim B_\mu$. We claim that 
\be
b_\mu=b_{-\mu},\;\;\forall \mu>0.
\label{eq: c3}
\ee
Indeed, if $\psi \in B_\mu\setminus\{0\}$ we deduce $Z\psi =0$ and $T\psi=\mu
\psi$.   The first equality implies that $\psi$ is covariant constant along the 
fibers of $N$ and thus $\psi$ is the pullback of some section $\hat{\psi}$ on $K^{-1}\otimes L \oplus L \ra \Sigma$. 
The second equality implies that $\hat{\psi}$ is a $\mu$-eigenvector of the    
${\bZ}_2$-graded, $L$-twisted,  Hodge-Dolbeault operator $\bar{\partial}+\bar{\partial}^*$ on $\Sigma$. The equality (\ref{eq: c3}) is now obvious.   Hence
\be
\eta(s)=\sum_{\mu>0}\frac{f_\mu-f_{-\mu}}{\mu^s}.
\label{eq: c4}
\ee
At this point, the dimensions $f_\mu$ can be described quite explicitly. More precisely  we have
\be
f_\mu \neq 0 \Rightarrow \mu \in {\bZ}\;\;{\rm and}\;\; f_\mu=h(K-L-\mu L_0)+h(L-\mu  L_0).
\label{eq: c5}
\ee
Before we prove the above equality we want to show its impact on the computation of $\eta(s)$. Using  it (\ref{eq: c4})  we deduce
\[
\eta(s)=\sum_{\mu>0}\frac{ (\,h(K-L-\mu L_0)-h(L+\mu L_0)\,)+(\,h(L-\mu L_0)-h(K-L+\mu   L_0)\,)}{\mu^s}
\]
\[
{\rm (Riemann-Roch)}=\sum_{\mu>0}\frac{ -(\deg L+\mu\deg L_0 +1-g)+(\deg L-\mu\deg L_0+1-g)}{\mu^s}
\]
\[
=\sum_{\mu>0}\frac{-2\mu}{\mu^s}=-2\ell\zeta(s-1)
\]
where $\zeta(s)$ is Riemann's zeta function. In particular, 
\[
\eta(0) =-2\ell \zeta(-1)
\]
while by [WW], $\zeta(-1)=-\frac{1}{12}$. This  agrees with (\ref{eq: adiaeta}).

\bigskip

\noindent{\bf Proof of} (\ref{eq: c5}).  Let $\psi \in F_\mu$. We decompose $\psi=\alpha\oplus \beta$ using the splitting (\ref{eq: c0}). Then
\be
\ii\nabla^B_\zeta \alpha =\mu \alpha, \;\; \bar{\partial}_B^*\alpha =0.
\label{eq: c6}
\ee
\be
-\ii\nabla_\zeta^B\beta=\mu\beta,\;\;\bar{\partial}_B\beta=0.
\label{eq: c7}
\ee
Denote by $f_\mu^-$ (resp. by $f_\mu^+$) the dimension of the space of solutions of (\ref{eq: c6})  (resp. (\ref{eq: c7})).  (The connections intervening in the above equations are connections on ${\can}^{-1}\otimes \pi^*L$ and $\pi^*L$ obtained by pullback from connections on  line bundles over the base.)   We will only show how to determine $f_\mu^+$ since the determination of $f_\mu^-$ is entirely similar.

Set for simplicity $\hat{B}=\pi^*B$, $\hat{B}_\mu^\pm= \hat{B}\mp \ii \mu\vfi$.   Note first that $\bar{\partial}_{\hat{B}}\beta =\bar{\partial}_{\hat{B}_\mu^\pm}$  since the transition $\hat{B}\ra \hat{B}_\mu^\pm$ does not alter the derivatives along horizontal directions. On the other hand, the equation $-\ii\nabla_\zeta^{\hat{B}}\beta=\mu \beta$ can be rewritten as
\[
-\ii\nabla_\zeta^{\hat{B_\mu^+}}\beta=0.
\]
Thus the equations (\ref{eq: c7}) are equivalent to
\be
-\ii\nabla_\zeta^{\hat{B_\mu^+}}\beta=0, \;\; \bar{\partial}_{B_\mu^+}\beta=0.
\label{eq: c8}
\ee
If  (\ref{eq: c8}) admits a nontrivial solution $\beta$ then $\beta$ must be $B_\mu^+$-covariant constant along the fibers. This implies that the pair $(\pi^*L, B_\mu^+)$ is the pullback of a pair (line bundle $L'$+connection
$B'$on $L'$) on the base $\Sigma$. The curvature of the connection $B'$ can be determined from
\[
F_{B'}=F_B+2\ell\mu\ii dvol_\Sigma.
\]
so that
\[
-\mu\ell +\deg L=\frac{\ii}{2\pi}\int_\Sigma F_{B'}\in {\bZ} 
\]
On the other hand, since $\pi^*L'\cong \pi^*L$ we deduce that $\mu\ell \cong 0$ (mod $\ell$) so that $\mu \in {\bZ}$.  In fact one can see that we have an isomorphism of {\em holomorphic line bundles}
\[
(L, B')\cong L\otimes L_0^{-\mu}.
\]
The second equation of (\ref{eq: c8}) implies that $\beta$ is a holomorphic section of $L\otimes L_0^{-\mu}$. Hence $f_\mu^+=h(L-\mu L_0)$. Similarly $f_\mu^-=h(K-L-\mu L_0)$ which concludes the proof of (\ref{eq: c5}).

\begin{remark}{\rm  A similar argument allows one to compute the  entire eta function of an adiabatic Dirac coupled with a flat connection of the type discussed in \S 2.4. In the notation of that section we have
\be
\eta(s)= -\ell \left( \zeta(s-1, k/|\ell|) +\zeta(1-s,1-k/|\ell|) \right)
\label{eq: rh}
\ee
where for any $a\in (0,1]$   we denoted by $\zeta(s,a)$ the Riemann-Hurwitz function
\[
\zeta(s,a) =\sum_{n=0}^\infty \frac{1}{(n+a)^s}
\]
In [WW] it is shown  that $\zeta(-1,a)= -B'_3(a)/6$ where $B'_3$ denotes the derivative of the Bernoulli polynomial $B_3(z)= z^3-\frac{3}{2}z^2 +\frac{1}{2}z$.   Substituting this in (\ref{eq: rh}) we reobtain the  main result of \S 2.4.}
\label{rem: rh}
\end{remark}

\section{Technical identities}
\setcounter{equation}{0}
We gather here  some technical results  used at various places in the main body of the paper.

Consider  a local orthonormal coframe  $\{\vfi, \vfi_1,\vfi_2\}$ of $T^*N$ as in \S 2.1. Set
\[
\ve =\frac{1}{\sqrt{2}}(\vfi_1+\ii\vfi_2)\;\;\;\bar{\ve}=\frac{1}{\sqrt{2}}(\vfi_1-\ii \vfi_2).
\]
Note that ${\ve}$ is a local section of ${\can}$. With respect to the  splitting 
\[
{\bS}_L={\can}^{-1/2}\otimes L\oplus {\can}^{1/2}\otimes L
\]
the Clifford multiplication has the block decomposition (see [N])
\[
{\bf c}(a\vfi +b\vfi_1+c\vfi_2)=\left[
\begin{array}{cc}
\ii a & (b+\ii c)\bar{\ve}\\
-(b-\ii c)\ve & -\ii a
\end{array}
\right]
\]
In particular, ${\bf c}(\vfi){\bf c}(\vfi_1){\bf c}(\vfi_2)=-1$ which agrees with the conventions described in Lemma 1.22 of [BC].

We compute easily
\[
{\bf c}({\ve})=\left[
\begin{array}{cc}
0 & \\
-\sqrt{2}{\ve} & 0
\end{array}
\right]\;\;{\bf c}(\bar{\ve})=\left[
\begin{array}{cc}
0 & \sqrt{2}{\ve} \\
0 & 0 
\end{array}
\right].
\]
If $\ii \da \in \ii\Omega^1(N)$  has the orthogonal decomposition
\[
\ii\da =\ii\da_0\vfi+\frac{1}{2}(\omega -\bar{\omega}), \;\;\omega \in
C^\infty({\can}),\;\;a_0\in C^\infty(N) 
\]
then
\be
{\bf c}(\ii \da)= \left[
\begin{array}{cc}
-\da_0 & -2^{-1/2}\bar{\omega} \\
-2^{-1/2}{\omega} & \da_0
\end{array}
\right].
\label{eq: d1}
\ee
The quadratic map  $q(\phi)$, viewed as an endomorphism of  ${\bS}_L$ has the block decomposition
\[
q(\phi)=\left[
\begin{array}{cc}
\frac{1}{2}(|\phi_-|^2-|\phi_+|^2) & \phi_-\bar{\phi}_+ \\
\bar{\phi}_-\phi_+ &  -\frac{1}{2}(|\phi_-|^2-|\phi_+|^2) 
\end{array}
\right]
\]
If  $\phi$ is such that $\phi_-=0$ then
\[
\dot{q}(\phi, \dps)\stackrel{def}{=}\frac{d}{dt}\!\mid_{t=0}q(\phi+t\dps)=\left[
\begin{array}{cc}
-{\re}\lan\phi_+, \dps_+\ran  & \dps_-\bar{\phi}_+ \\
\dbps_-\phi_+ & {\re}\lan \phi_+, \dps_+ \ran
\end{array}
\right].
\]
Using (\ref{eq: d1}) we obtain  a description of $\dot{q}(\phi, \dps)$ as  a {\em purely imaginary
1-form}
\be
\dot{q}(\phi, \dps) = \ii{\re}\lan \phi_+,\dps_+ \ran\vfi -2^{-1/2}(\dps_-\bar{\phi}_+ -\dbps_-\phi_+).
\label{eq: d2}
\ee
Let $\phi$ as above. Given  $\Xi=\dps \oplus \ii \da \oplus \ii f$ where $\ii \da =\frac{1}{2}(\omega-\bar{\omega})$, $\omega\in C^\infty({\can})$  then using (\ref{eq: d1}) and (\ref{eq: d2}) we deduce
\be
{\cal P}_\phi \Xi =\left[
\begin{array}{c}
{\bc}(\ida)\phi -\ii f \phi \\
\dta(\phi, \dps) \\
\ii \im \lan \phi, \dps\ran
\end{array}
\right]=\left[
\begin{array}{c}
(-\bar{\omega}\phi_+ )\oplus (-\ii f\phi_+) \\
= \ii{\re}\lan \phi_+,\dps_+ \ran\vfi +2^{-1/2}(\dps_-\bar{\phi}_+ -\dbps_-\phi_+)\\
\ii {\im}\lan \phi_+, \dps_+\ran
\end{array}
\right]
\label{eq: d3}
\ee
We now deduce easily that
\be
\lan {\cal P}\Xi, \Xi\ran = 2^{1/2}f{\im}\lan \phi_+, \dps_+\ran - {\re}(\dbps_-\phi_+ \bar{\omega}).
\label{eq: d4}
\ee
The term $\dot{q}(\phi, \dps)$  has nice  divergence properties. More precisely we have the following result.

\begin{lemma}{\rm   Consider a  $spin^c$ structure $\si$ on an oriented,
 Riemannian $3$-manifold $(M, g)$.  Fix a  connection $A$ on $\det \si$  and denote by
 ${\dir}_A$ the  Dirac operator on ${\bS}_\sigma$  induced by the
 {\em Levi-Civita connection} coupled with $A$. Then for every $\psi \in C^\infty ({\bS}_\si)$ we have} 
 \[
 d^* q(\psi)=-\ii \im \lan \psi, {\dir}_A\psi\ran.
 \]
 \label{lemma: divergence}
 \end{lemma}

 \noindent{\bf Proof} \hspace{.3cm}   Fix an arbitrary point $p_0\in M$,  choose  normal
 coordinates $(x^1, x^2, x^3)$ near $p_0$ and set $e^i=dx^i$. Note that   {\em
 at $p_0$} we have  $d^*e^i=0$ for  all $i$. In [N] we showed that, viewed  as a
 1-form, $q(\psi)$ has the local description
 \[
 q(\psi)=\frac{1}{2}\sum_i\lan \psi, \bc(e^i) \psi\ran e^i.
 \]
 {\em At $p_0$ } we have
 \[
 2d^* q(\psi)=-\sum_i\pai(\lan\psi, \bc (e^i) \psi \ran) e^i
 \]
 \[ 
 =-\sum_i\lan\nabla^A_i\psi, \psi \ran  -\sum_i \lan \psi, \bc(e^i) \nabla^A_i \psi
 \ran\;\;\;({\rm since}\; \nabla_ie^i=0\;{\rm at}\; p_0) 
 \]
 \[
 =\sum_i \overline{\lan \psi,\bc(e^i)\nabla^A_i\psi\ran}-\sum_i\lan \psi,
 \bc(e^i)\nabla_i^A \psi\ran =-2{\bf i}\im \lan \psi, {\dir}_A\psi \ran .
 \]
 Since $p_0$ is arbitrary this proves the lemma. $\Box$
 
 \bigskip
 
 On our circle bundle $N$ we have  ${\dir}_A= D_A+\lambda_r/2$ so that
 \be
 d^* q(\phi)=-\ii \im \lan\,\phi, (D_A+\lambda_r/2)\phi \,\ran =-\ii \im \lan \phi ,
 D_A \phi \ran.
 \label{eq: divergence}
 \ee
Suppose now $\phi\in\ker D_A$.  We derivate (\ref{eq: divergence}) along $\dps$  and we  get
\be
d^*\dot{q}(\phi, \dps)=-\im \lan\phi, D_A\dps\ran.
\label{eq: d5}
\ee
This identity  plays  an important role in the proof of the following  result.

\begin{lemma}{\rm  Consider   an  irreducible solution $\co=(\phi, A)$ of the  Seiberg-Witten equation on $N$. We assume for simplicity $\phi_-=0$. For each $w \geq 0$ we have an operator ${\cal O}_w={\cal O}_w(\co)$ as in \S 3.4.   (Recall that ${\cal O}_0={\hhes}_1$.) Then for all $w\geq 0$ we have}
\[
\ker {\hhes}_1=\ker {\cal O}_w.
\]
\label{lemma: d1}
\end{lemma}

\noindent{\bf Proof}\hspace{.3cm}   Note that if $\Xi=\dps\oplus \ii\da \oplus
\ii f \in \ker {\cal O}_w$ is such that  $f\equiv 0$ then the
  definition of ${\cal O}_w$ implies immediately that $\Xi\in \ker {\hhes}_1$. 
  Conversely, any $\Xi\in {\hhes}_0$ has vanishing third coordinate. Hence it suffices
  to show that the third component of any $\Xi\in \ker {\cal O}_w$ vanishes.  

Let $\Xi=\dps\oplus \ii \da \oplus \ii f\in \ker {\cal O}_w$. This means
\be
\left\{
\begin{array}{clclcl}
D_A\dps  &                               &+ &  {\bc}(\ida)\phi-\ii f\phi &=&0  \\
                 &-\ii\ast  d\da +\ii df & +&\dta(\phi, \dps) &=&0\\
                  &\ii d^*\da -2w\ii f& + & \ii \im \lan \phi, \dps\ran &=&0 
\end{array}
\right.
\label{eq: d7}
\ee
Take the inner product of the second equation with $\ii df$. After an integration by parts we get
\[
\int_N |df|^2 dv_N -\int_N f\cdot d^*\ii \dot{q}(\phi, \dps dv_N=0.
\]
Using (\ref{eq: divergence}) and the first equation in (\ref{eq: d7})  and we get
\[
\int_N (|df|^2 +|f|^2\cdot|\phi|^2) dv_N=0.
\]
This shows $f\equiv 0$ and completes the proof of the lemma. $\Box$

\bigskip

The above lemma has the following important consequence.

\begin{corollary}
\[
SF({\hhes}_1\ra {\cal O}_w)=0.
\]
\label{cor: d1}
\end{corollary}

\vspace{5cm}

{\bf Current address}: Dept.of Math., McMaster University, Hamilton, Ontario,  L8S 4K1, Canada; nicolaes@icarus.math.mcmaster.ca

\bigskip

 {\bf Address beginning Sept. 1998}: Dept. of Math., University of Notre Dame, Notre Dame, IN 46556-5683

\end{document}